\documentclass[a4paper,11pt,english]{article}
\usepackage{babel}
\usepackage[latin1]{inputenc}
\author{Andrea Tellini \\ \small{Centre d'Analyse et de Math\'ematique Sociales} \\ \small{Ecole des Hautes Etudes en Sciences Sociales} \\ \small{190-198 Avenue de France -- 75013 Paris, France} \\ \texttt{andrea.tellini@ehess.fr}}
\title{\textbf{Propagation speed in a strip bounded by a line with different diffusion}\thanks{This work was mainly supported by the Spanish Ministry of Economy and Competitiveness through grants BES-2010-039030 and EEBB-I-13-05962. Partial founding was also provided by European Union's ERC Grant Agreement n. 321186 - ReaDi - Reaction-Diffusion Equations, Propagation and Modelling and Research Project
\lq\lq Stabilità asintotica di fronti per equazioni paraboliche\rq\rq of the University of Padova (2011) coordinated by Luca Rossi.}}
\date{\today}
\usepackage{amsmath,amsfonts,amssymb}
\usepackage{yhmath,mathrsfs,yfonts,arcs}
\usepackage{graphicx}

\usepackage[bookmarks=false]{hyperref}
\usepackage{amsthm}
\usepackage{cleveref}

\newtheorem{theorem}{Theorem}[section]
\newtheorem{remark}[theorem]{Remark}
\newtheorem{proposition}[theorem]{Proposition}

\newtheorem{definition}[theorem]{Definition}

\DeclareMathOperator{\KPP}{KPP}
\DeclareMathOperator{\DR}{DR}
\DeclareMathOperator{\Int}{int}

\newcommand{\field}[1] {\mathbb{#1}}
\newcommand{\N}{\field{N}}

\newcommand{\R}{\field{R}}
\newcommand{\C}{\field{C}}

\newcommand{\cKPP}{c_{\KPP}}
\newcommand{\cint}{c_{\Int}}

\def\a{\alpha}
\def\b{\beta}
\def\e{\varepsilon}
\def\D{\Delta}
\def\d{\delta}
\def\g{\gamma}

\def\m{\mu}
\def\n{\nu}

\def\o{\omega}
\def\O{\Omega}
\def\p{\partial}
\def\r{\rho}

\def\S{\Sigma}
\def\s{\sigma}
\def\t{\theta}

\def\ov{\overline}
\def\un{\underline}
\def\ua{\uparrow}
\def\da{\downarrow}

\newcommand{\mc}{\mathcal}

\begin{document}
\maketitle

\vspace{-1.0cm}
\begin{center}
\noindent
\emph{In memory of Giuliano Bardi (1948--2012),\\ the one who taught me what a derivative is.}
\end{center}

\begin{abstract}
In this paper we consider a model for the diffusion of a population in a strip-shaped field, where the growth of the species is governed by a Fisher-KPP equation and which is bounded on one side by a road where the species can have a different diffusion coefficient. Dirichlet homogeneous boundary conditions are imposed on the other side of the strip. We prove the existence of an asymptotic speed of propagation which is greater than the one of the case without road and study its behavior for small and large diffusions on the road. Finally we prove that, when the width of the strip goes to infinity, the asymptotic speed of propagation approaches the one of an half-plane bounded by a road, case that has been recently studied in \cite{BRR1, BRR2}.
\end{abstract}



\setcounter{equation}{0}
\setcounter{figure}{0}
\section{Introduction}
\label{sec1}
Recently, in \cite{BRR1}, the system
\begin{equation}
\label{11}
  \left\{ \begin{array}{lll}
  \!\!\!u_t(x,t)\!-\!Du_{xx}(x,t)\!=\!\nu v(x,L,t)\!-\!\m u(x,t) & \!\!\text{for } x\in\R, & \!\! t>0 \\
  \!\!\!v_t(x,y,t)-d\D v(x,y,t)=f(v) & \!\!\text{for } \!(x,y)\!\in\!\R\!\times\!(-\infty,L), & \!\! t>0  \\ 
  \!\!\!dv_y(x,L,t)=\m u(x,t)-\nu v(x,L,t) & \!\!\text{for } x\in\R, & \!\! t>0.
  \end{array} \right.
\end{equation}
was introduced to model the evolution of the species $v(x,y,t)$ in a \emph{field} $\R\times(-\infty,L)$ which is bounded at the level $y=L$ by a \emph{road} where part of the same species, $u(x,t)$, diffuses with coefficient $D>0$, which in principle may be different from the diffusion coefficient in the field $d>0$. A reaction of Fisher-KPP type takes place in the field, i.e. $f\in\mc{C}^1([0,+\infty))$ satisfies
\begin{equation}
\label{12}
f(0)=0=f(1), \qquad 0<f(s)<f'(0)s \;\; \text{ in } (0,1), \qquad f<0 \text{ in } (1,+\infty).
\end{equation}
On the contrary, no reaction occurs on the road, where the density of the species varies only because a fraction $\mu>0$ of the population jumps from the road to the field while a fraction $\nu>0$ of the population jumps from the field to the road.

This model was motivated by empirical observations of wolves moving along seismic lines in Canada (see \cite{McK}) or insects like the \emph{Aedes albopictus} (tiger mosquito) spreading in the United States along highways (see \cite{MM}). Another example of this phenomenon is the diffusion of diseases along commercial and transport networks (see \cite{TRH} and the references therein).

In \cite{BRR1}, the authors established the existence of an asymptotic speed of propagation (see Definition \ref{deasp}) of the solution of \eqref{11}, starting from a continuous, nonnegative, compactly supported initial datum $(u_0, v_0)\neq(0, 0)$, towards the unique steady state of the problem, as well as some qualitative properties of it. Denoting such speed by $c^*_{\infty}$, they showed that, if $D\leq 2d$, then $c^*_{\infty}=\cKPP$, where
\begin{equation}
\label{13}
\cKPP=2\sqrt{df'(0)}
\end{equation}
is the asymptotic speed of propagation of the classical Fisher-KPP equation
\begin{equation*}
v_t(x,y,t)-d\D v(x,y,t)=f(v(x,y,t))
\end{equation*}
in the half-plane (see \cite{KPP,AW}), while, if $D> 2d$, then $c^*_{\infty}>\cKPP$. This means that a large diffusion on the road speeds up the propagation of the population in the field. Moreover the authors showed that the spreading velocity increases to infinity as the diffusivity on the line grows to infinity. In \cite{BRR2} they also studied the influence that  a drift term and a Fisher-KPP reaction also on the road have on the asymptotic speed of propagation.

In this work we investigate the effect of the road on the propagation in a field which is no longer a half-plane but a strip $\O=\R\times(0,L)$. On the other part of the boundary of the field we impose homogeneous Dirichlet boundary conditions, modeling in this way an unfavorable region at level $y=0$. The system we consider is therefore
\begin{equation}
\label{14}
  \left\{ \begin{array}{lll}
  u_t(x,t)-Du_{xx}(x,t)=\nu v(x,L,t)-\m u(x,t) & \text{for } x\in\R, & t>0 \\
  v_t(x,y,t)-d\D v(x,y,t)=f(v) & \text{for } (x,y)\in\O, & t>0  \\ 
  dv_y(x,L,t)=\m u(x,t)-\nu v(x,L,t) & \text{for } x\in\R, & t>0 \\
  v(x,0,t)=0 & \text{for } x\in\R, & t>0, \end{array}\right.
\end{equation}
where $D,d,\mu,\nu,L$ are positive constants. Ascertaining the long time behavior of the solutions of the Cauchy problem associated to \eqref{14} will be the first step in the study of the problem. The result is the following
\begin{theorem}
\label{thltb}
Let $(u,v)$ denote the solution of \eqref{14} starting from a nonnegative, not equal to $(0,0)$, bounded and continuous initial datum $(u_0, v_0)$. If
\begin{equation}
\label{15}
\frac{f'(0)}{d}\leq\left(\frac{\pi}{2L}\right)^2,
\end{equation}
then 
\begin{equation}
\label{16}
\lim_{t\to+\infty}(u(x,t),v(x,y,t))=(0,0)
\end{equation}
locally uniformly in $\overline{\O}$, while if
\begin{equation}
\label{17}
\frac{f'(0)}{d}>\left(\frac{\pi}{2L}\right)^2
\end{equation}
and, moreover, 
\begin{equation}
\label{18}
\frac{f(s)}{s} \text{ is nonincreasing},
\end{equation}
then
\begin{equation}
\label{19}
\lim_{t\to+\infty}(u(x,t),v(x,y,t))=\left(\frac{\nu}{\mu}V(L),V(y)\right)
\end{equation}
locally uniformly in $\overline{\O}$, where $V(y)$ is the unique solution of
\begin{equation}
\label{110}
  \left\{ \begin{array}{l}
  -dV''(y)=f(V(y)) \quad \text{for } y\in(0,L)\\
  V'(L)=0, \quad V(0)=0. \end{array}\right.
\end{equation}
\end{theorem}
A remarkable difference with respect to \eqref{11} is that here the width of the strip $L$ plays a role in the existence of positive steady states. Moreover, condition \eqref{18} was not necessary to guarantee uniqueness in \cite{BRR1}, while here it is (see Remark \ref{re33} below). From a biological point of view, Theorem \ref{thltb} says that, if the strip is not sufficiently large,  the influence of the unfavorable region $y=0$ drives the species to extinction. On the contrary, the species will persist if the strip is sufficiently large. In this latter case, a natural question is to study more deeply how the convergence to the steady state occurs. To this end we consider the following concept
\begin{definition}
\label{deasp}
We say that $c^*\in\R_+$ is an \emph{asymptotic speed of propagation} (in the $x-$direction) for \eqref{14} if, denoting by $(u, v)$ the solution of \eqref{14} with
a continuous, nonnegative, compactly supported initial datum $(u_0, v_0)\neq(0, 0)$, we have
\begin{itemize}
\item for all $c > c^*$,
\begin{equation}
\label{111}
\lim_{t\to +\infty} \sup_{\substack{|x|\geq ct \\ y\in[0,L]}} |(u(x, t), v(x, y, t))| = 0,
\end{equation}
\item for all $0<c < c^*$,
\begin{equation}
\label{112}
\lim_{t\to +\infty}\sup_{\substack{|x|\leq ct \\ y\in[0,L]}} \left|(u(x, t), v(x, y, t)) - \left(\frac{\nu}{\mu}V(L), V(y)\right)\right|=0,
\end{equation}
where $V(y)$ is the unique solution of \eqref{110}. 
\end{itemize}
\end{definition}
\noindent
In this sense, the main result of the paper is the following
\begin{theorem}
\label{thmain}
Problem \eqref{14} admits an asymptotic speed of propagation, denoted by $c^*=c^*_L(D,d,\mu,\nu)$, such that:
\begin{itemize}
\item[(i)] $D\mapsto c^*(D)$ is increasing,
\item[(ii)] the following limits exist and are positive real numbers
\begin{equation*}
\lim_{D\to 0}c^*(D)=\ell_0,
\qquad
\lim_{D\to+\infty}\frac{c^*(D)}{\sqrt{D}}=\ell_\infty,
\end{equation*}
\item[(iii)] for fixed $D,d,\mu,\nu$ we have
\begin{equation*} 
\lim_{L\to+\infty}c^*_{L}(D,d,\mu,\nu)=c^*_{\infty}(D,d,\mu,\nu),
\end{equation*}
where $c^*_{\infty}(D,d,\mu,\nu)$ is the asymptotic speed of propagation of Problem \eqref{11}.
\end{itemize}
\end{theorem}

The paper is organized like follows: in Section \ref{sec2} we recall some tools from \cite{BRR1,BRR2} which will be indispensable throughout the rest of the paper; Section \ref{sec3} provides the proof of Theorem \ref{thltb}; in Sections \ref{sec4} and \ref{sec5} we construct $c^*$ and derive some properties that will allow us, in Section \ref{sec6}, to show that it satisfies Definition \ref{deasp} and relation $(i)$ of Theorem \ref{thmain}. In Section \ref{sec7} we prove relation $(ii)$ of Theorem \ref{thmain} and finally in Section \ref{sec8} we study the influence of the road on the asymptotic speed of propagation, comparing \eqref{14} with the case in which no road is present, and give the proof of Theorem \ref{thmain}$(iii)$.

\setcounter{equation}{0}
\setcounter{figure}{0}
\section{Preliminary results}
\label{sec2}
In this section we present some fundamental results that are contained in or follow easily from \cite{BRR1,BRR2}.
The existence of a solution for the Cauchy problem associated to \eqref{14} with a continuous initial datum $(u_0,v_0)$ follows from an easy modification of \cite[Appendix A]{BRR1} and uniqueness follows from a comparison principle which will be diffusely used throughout this paper and whose proof can be easily adapted from \cite[Proposition 3.2]{BRR1}. Before stating it we point out that, as usual, by a \emph{supersolution} (resp. \emph{subsolution}) of \eqref{14} we mean a pair $(u,v)$ satisfying System \eqref{14} with $\geq$ (resp. $\leq$) instead of $=$.

\begin{proposition}
\label{prcp1}
Let $(\un{u}, \un{v})$ and $(\ov{u},\ov{v})$ be, respectively, a subsolution bounded from
above and a supersolution bounded from below of \eqref{14} satisfying $\un{u}\leq\ov{u}$ and $\un{v}\leq\ov{v}$
at $t = 0$. Then, either $\un{u}<\ov{u}$ and $\un{v}<\ov{v}$ for all $t$, or there exists $T > 0$ such that
$(\un{u}, \un{v})=(\ov{u},\ov{v})$ for $t\leq T$.
\end{proposition}

We also need the following comparison principle regarding an extended class of generalized subsolutions and which is a particular instance of \cite[Proposition 3.3]{BRR1}.

\begin{proposition}
\label{prcp2}
Let $E\subset \R\times\R_+$ and $F\subset \O\times\R_+$ be two open sets,
let $(u_1, v_1)$ be a subsolution of \eqref{14} bounded from above and satisfying
\begin{equation*}
u_1 = 0 \text{ on } (\p E) \cap (\R\times\R_+), \quad v_1= 0 \text{ on } (\p F) \cap (\O\times\R_+),
\end{equation*}
and consider
\begin{equation*}
\un{u} :=\begin{cases}
\max\{u_1, 0\} & \text{in $\ov{E}$} \\
0 & \text{otherwise,}
\end{cases} 
\qquad
\un{v} :=\begin{cases}
\max\{v_1, 0\} & \text{in $\ov{F}$} \\
0 & \text{otherwise.}
\end{cases}
\end{equation*}
If they satisfy
\begin{equation}
\label{21}
\begin{split}
\un{v}(x, L, t) &\geq v_1(x, L, t) \quad \text{ for all } (x,t) \text{ such that } \un{u}(x, t) > 0, \\
\un{u}(x, t) &\geq u_1(x, t) \qquad \text{ for all } (x,t) \text{ such that } \un{v}(x, L, t) > 0,
\end{split}
\end{equation}
then, for any supersolution $(\ov{u},\ov{v})$ of \eqref{14} bounded from below and such that $\un{u}\leq\ov{u}$ and
$\un{v}\leq\ov{v}$ at $t = 0$, we have $\un{u}\leq\ov{u}$ and $\un{v}\leq\ov{v}$ for all $t > 0$.
\end{proposition}

\begin{remark}
\label{recp2}
The same result of Proposition \ref{prcp2} holds for problems like \eqref{14} with an additional drift term in the differential operator.
\end{remark}

As a consequence of the previous analysis, we will consider continuous nonnegative initial data throughout the rest of this work, since we are interested in nonnegative solutions of \eqref{14}.

\setcounter{equation}{0}
\setcounter{figure}{0}
\section{Liouville-type result and long time behavior}
\label{sec3}
In order to determine the long time behavior of the solutions of \eqref{14} we need to study the solutions of the elliptic system associated to it, precisely
\begin{equation}
\label{31}
  \left\{ \begin{array}{ll}
  -D U_{xx}(x)=\nu V(x,L)-\m U(x) & \text{for } x\in\R, \\
  -d\D V(x,y)=f(V) & \text{for } (x,y)\in\O,   \\ 
  d V_y(x,L)=\m U(x)-\nu V(x,L) & \text{for } x\in\R, \\
  V(x,0)=0 & \text{for } x\in\R. \end{array}\right.
\end{equation}
Actually, Propositions \ref{pr31} and \ref{pr34} below suggest that we have to focus on solutions of \eqref{31} which are $x-$independent. They are of the form $(U,V(y))$, where $V$ satisfies \eqref{110} and, thanks to the first equation of \eqref{31}, $U=\frac{\nu}{\m} V(L)$.
The first result regarding the long time behavior is the following
\begin{proposition}
\label{pr31}
Let $(u,v)$ be the solution of \eqref{14} starting with a nonnegative, bounded initial datum $(u_0,v_0)$. Then, there exists a nonnegative, bounded solution $V_1$ of \eqref{110} such that
\begin{equation*}
\limsup_{t\to+\infty}u(x,t)\leq U_1, \qquad \limsup_{t\to+\infty}v(x,y,t)\leq V_1(y)
\end{equation*}
locally uniformly in $\overline{\O}$, where
\begin{equation*}
U_1=\frac{\nu}{\m} V_1(L).
\end{equation*}
\begin{proof}
Observe preliminarily that, if we define, for $(x,y)\in\overline{\O}$,
\begin{equation*}
\overline{v}(x,y)=\max\left\{1,\left\|v_0\right\|_{\infty},\frac{\mu}{\nu}\left\|u_0\right\|_{\infty}\right\}, \qquad  \overline{u}(x)=\frac{\nu}{\mu}\overline{v},
\end{equation*}
then $(\overline{u},\overline{v})$ is a strict supersolution of \eqref{31} which is larger than $(u_0,v_0)$. Therefore, by Proposition \ref{prcp1}, we have that the solution of \eqref{14} with $(\overline{u},\overline{v})$ as initial datum is decreasing and, thanks to parabolic estimates, converges locally uniformly in $\overline{\O}$ to a nonnegative stationary solution $(U_1,V_1)$ of \eqref{14}, i.e. a solution of \eqref{31}. Proposition \ref{prcp1} also gives
\begin{equation*}
\limsup_{t\to+\infty} u(x,t)\leq U_1(x), \qquad \limsup_{t\to+\infty} v(x,y,t)\leq V_1(x,y).
\end{equation*}
From the invariance of Problem \eqref{14} in the $x$ direction and the uniqueness of the associated Cauchy problem, translations in $x$ of a solution of \eqref{14} with a certain initial datum coincide with the solution of \eqref{14} starting from the translated initial datum. Since $(\ov{u},\ov{v})$ is $x-$independent, the $x-$invariance of $(U_1,V_1)$ follows. 
\end{proof}
\end{proposition}

Obviously, \eqref{110} admits the trivial solution $V=0$. In the following proposition we will show that \eqref{17} is a necessary and sufficient condition for \eqref{110} to possess positive bounded solutions.

\begin{proposition}
\label{pr32}
Problem \eqref{110} admits positive bounded solutions if and only if \eqref{17} holds.
Moreover, if we assume \eqref{17} and \eqref{18}, then Problem \eqref{110} admits a \emph{unique} positive bounded solution.

\begin{proof}
We begin with the necessity of \eqref{17}. Suppose it does not hold and Problem \eqref{110} admits a positive solution $v$. Then, multiplying the differential equation of \eqref{110} by $\sin(\frac{\pi}{2L}y)$ and integrating by parts in $(0,L)$, we get
\begin{multline*}
\int_0^Lf(v(y))\sin\left(\frac{\pi}{2L}y\right)\,dy=d\left(\frac{\pi}{2L}\right)^2\int_0^Lv(y)\sin\left(\frac{\pi}{2L}y\right)\,dy\geq \\
\geq f'(0)\int_0^Lv(y)\sin\left(\frac{\pi}{2L}y\right)\,dy>\int_0^Lf(v(y))\sin\left(\frac{\pi}{2L}y\right)\,dy,
\end{multline*}
where, for the last relation, we have used the second assumption in \eqref{12}. We have reached a contradiction and therefore no positive solution can exist.

Now we pass to the sufficiency. First of all we show that any positive solution of \eqref{110} must satisfy $v(y)<1$ in $[0,L]$. Indeed if there was a point $y_0\in(0,L]$ where $v(y_0)=1$, either it would be a maximum of $v$ in a (relative to $(0,L]$) neighborhood  of $y_0$, but in such a case $v\equiv 1$ by the uniqueness of the Cauchy problem $-v''=f(v)$ with conditions $v(y_0)=1, v'(y_0)=0$, or there would be $y_1\in(0,L]$ with $v(y_1)=\max v>1$. But in this case, $-dv''(y_1)=f(v(y_1))<0$, which is impossible for a maximum.

As a consequence of this result and \eqref{12}, we have that $v''<0$ in $(0,L)$ and therefore $v'$ is decreasing. This means that $v'$ is positive in $(0,L)$, since $v'(L)=0$, i.e. $v$ is increasing in $(0,L)$. By multiplying the differential equation of \eqref{110} by $v'$ and integrating in $(y,L)$, with $0<y<L$, we get
\begin{equation*}
d\frac{v'(y)^2}{2}=\int_{v(y)}^{v(L)}f(s)\,ds
\end{equation*}
and, recalling that $v'>0$, we have that any solution of \eqref{110} must satisfy
\begin{equation*}
L=\int_0^L\frac{v'(y)\,dy}{\sqrt{\frac{2}{d}\int_{v(y)}^{v(L)}f(s)\,ds}}=\int_0^1\frac{d\xi}{\sqrt{\frac{2}{d}\int_{\xi}^{1}\frac{f(v(L)\eta)}{v(L)\eta}\eta\,d\eta}}
\end{equation*}
with $0<v(L)<1$. Therefore, if we define the function
\begin{equation}
\label{33}
M(\r):=\int_0^1\frac{d\xi}{\sqrt{\frac{2}{d}\int_{\xi}^{1}\frac{f(\r\eta)}{\r\eta}\eta\,d\eta}},
\end{equation}
which is continuous in $(0,1)$ and measures the length of the interval necessary for a solution of
\begin{equation*}
  \left\{ \begin{array}{l}
  -dV''(y)=f(V(y)) \quad \text{for } y\in(0,y_0)\\
  V(0)=0, \quad V'(y_0)=0 \end{array}\right.
\end{equation*}
to attain its maximum value $\r$ (at $y_0$), the uniqueness of the Cauchy problem associated to the ordinary differential equation of \eqref{110}, we have that any solution of $M(\r)=L$ provides a solution of \eqref{110}. The function $M$ satisfies
\begin{equation*}
\lim_{\r\ua 1}M(\r)=+\infty
\end{equation*}
since a maximum equal to $1$ cannot be attained in a finite interval, as seen before. Moreover, thanks to \eqref{17},
\begin{equation*}
\lim_{\r\da 0}M(\r)=\sqrt{\frac{d}{f'(0)}}\int_0^1\frac{d\xi}{\sqrt{1-\xi^2}}=\sqrt{\frac{d}{f'(0)}}\frac{\pi}{2}<L
\end{equation*}
and, therefore, there exists $\bar{\r}\in(0,1)$ such that $M(\bar{\r})=L$, which provides us with a solution of \eqref{110}.

As far as uniqueness, it is easily seen that, under hypothesis \eqref{18}, the function $M$ is increasing and therefore there exists a unique value of $\rho$ for which $M(\rho)=L$.
\end{proof}
\end{proposition}

\begin{remark}
\label{re33}
If condition \eqref{18} does not hold, Problem \eqref{110} may exhibit more than one solution. Consider indeed $f(s)=s(-6s^3+9s^2-4s+1)$, which satisfies \eqref{12} but not \eqref{18}. With this choice, the function $M$ defined in \eqref{33}, which is $\mc{C}^1[0,1)$, satisfies
\begin{equation*}
M'(\r)=\sqrt{\frac{15\, d}{2}}\int_0^1h(\r,\xi)\,d\xi
\end{equation*}
where
\begin{equation}
\label{34}
h(\r,\xi):=\frac{216\r^2(1-\xi^5)-270\r(1-\xi^4)+80(1-\xi^3)}{\left(-72\r^3(1-\xi^5)+135\r^2(1-\xi^4)-80\r(1-\xi^3)+30(1-\xi^2)\right)^{3/2}}
\end{equation}
and, as a consequence, $M'(0)>0$, since $h(0,\xi)>0$ for every $\xi\in(0,1)$. On the other hand, for $\r=1/2$ the numerator in \eqref{34} reduces to
\begin{equation*}
h_1(\xi)=-54\xi^5+135\xi^4-80\xi^3-1,
\end{equation*}
which satisfies $h_1(0)=-1$, is decreasing for $\xi\in(0,2/3)$ and increasing for $\xi>2/3$. Since $h_1(1)=0$, this implies that $h(1/2,\xi)<0$ for all $\xi\in(0,1)$ and, therefore, $M'(1/2)<0$. Recalling that $M(\r)\to+\infty$ as $\r\ua 1$, the previous analysis entails that $L$ can be chosen in such a way that Problem \eqref{110} possesses at least $3$ solutions.
\end{remark}

The last result we need in order to prove Theorem \ref{thltb} is the following
\begin{proposition}
\label{pr34}
Assume \eqref{17} and let $(u,v)$ be the solution of \eqref{14} starting with a nonnegative, not equal to $(0,0)$, bounded initial datum. Then, there exists a positive bounded solution $V_2$ of \eqref{110} such that
\begin{equation*}
U_2\leq \liminf_{t\to+\infty}u(x,t), \qquad V_2\leq \liminf_{t\to+\infty}v(x,y,t) 
\end{equation*}
locally uniformly in $(x,y)\in\overline{\O}$, where
\begin{equation*}
U_2=\frac{\nu}{\mu}V_2(L).
\end{equation*}

\begin{proof}
Consider
\begin{equation}
\label{35}
(\un{u},\un{v}):=\left\{\begin{array}{ll}
\cos(\o x)\left(1,\frac{\mu\sin(\b y)}{d\b\cos(\b L)+\nu\sin(\b L)}\right) &\text{if $(x,y)\in\left(-\frac{\pi}{2\o},\frac{\pi}{2\o}\right)\times\left[0,L\right]$} \\
(0,0) &\text{otherwise}. 
\end{array}\right.
\end{equation}
We now show that $\b$ and $\o$ can be chosen so that $(\un{u},\un{v})$ satisfy
\begin{equation}
\label{36}
  \left\{ \begin{array}{ll}
  -D\un{u}_{xx}(x)\leq\nu \un{v}(x,L)-\m \un{u}(x) & \text{for } x\in \left(-\frac{\pi}{2\o},\frac{\pi}{2\o}\right) \\
  -d\D \un{v}(x,y)\leq(f'(0)-\d)\un{v} & \text{for } (x,y)\in\left(-\frac{\pi}{2\o},\frac{\pi}{2\o}\right)\times\left(0,L\right)  \\ 
  d\un{v}_y(x,L)=\m \un{u}(x)-\nu \un{v}(x,L) & \text{for } x\in \left(-\frac{\pi}{2\o},\frac{\pi}{2\o}\right) \\
  \un{v}(x,0)=0 & \text{for } x\in \left(-\frac{\pi}{2\o},\frac{\pi}{2\o}\right), \end{array}\right.
\end{equation}
for $0<\d<f'(0)$ and therefore, by the second relation of \eqref{12}, there exists $\e_0$ such that, for all $0<\e<\e_0$, $\e(\un{u},\un{v})$ is a strict generalized subsolution of \eqref{31} to which Proposition \ref{prcp2} can be applied. Observe that, with choice \eqref{35}, the last two equations of \eqref{36} are satisfied and the first two inequalities reduce to
\begin{equation*}
D\o^2\leq \frac{-\mu d\b\cos(\b L)}{d\b\cos(\b L)+\nu\sin(\b L)}, \qquad
d\o^2+d\b^2\leq f'(0)-\d.
\end{equation*}
Now, thanks to \eqref{17}, it is possible to fix $\d$ in such a way that
\begin{equation*}
d\left(\frac{\pi}{2L}\right)^2<f'(0)-\d
\end{equation*}
and take $\frac{\pi}{2L}<\b<\frac{\pi}{L}$ in a neighborhood of $\frac{\pi}{2L}$ (we denote $\b\sim\frac{\pi}{2L}$), so that
\begin{equation*}
m:=\min\left\{\frac{f'(0)-\d}{d}-\b^2,\frac{-\mu d\b\cos(\b L)}{D\left(d\b\cos(\b L)+\nu\sin(\b L)\right)}\right\}>0.
\end{equation*}
As a consequence, if $\o^2\leq m$, $(\un{u},\un{v})$ satisfies \eqref{36}.
Moreover, reducing $\e$ if necessary, we can assume that $\e(\un{u},\un{v})<(u(x,1),v(x,y,1))$, because, thanks to Proposition \ref{prcp1} and the Hopf lemma, 
we have that
\begin{equation*}
u(x,1)>0, \quad v_y(x,0,1)>0 \quad \text{ and }  \quad v(x,y,1)>0 \quad \text{ for all } (x,y)\in\O,
\end{equation*}
and, in addition, $v(x,L,1)>0$ for every $x\in\R$, since if there was $x_0$ such that $v(x_0,L,1)=0$, then from the third equation in \eqref{14} we would have
\begin{equation*}
dv_y(x_0,L,1)=\mu u(x_0,1)>0,
\end{equation*}
which is impossible, since, again by Proposition \ref{prcp1}, $v(x,y,1)\geq 0$ in $\ov{\O}$.

By Proposition \ref{prcp2}, the solution of \eqref{14} with $\e(\un{u},\un{v})$ as initial datum, converges, increasingly, to a stationary solution $(U_2,V_2)$ of \eqref{14} locally uniformly in $\overline{\O}$ and moreover
\begin{equation*}
U_2\leq\liminf_{t\to+\infty}u(x,t) \quad \text{ and } \quad V_2\leq\liminf_{t\to+\infty}v(x,y,t).
\end{equation*}
As before we have, for all $(x,y)\in\O$
\begin{equation*}
U_2(x)>0, \quad V_{2,y}(x,0)>0, \quad V_2(x,y)>0 \quad \text{ and }  \quad V_2(x,L)>0
\end{equation*}
and, since $\e(\un{u},\un{v})$ is continuous and compactly supported and, thanks to the above-mentioned monotonicity, it does not touch $(U_2,V_2)$, we have that there exists $k>0$, $k\sim 0$, such that $\e(\un{u}(x-h),\un{v}(x-h,y))$, which still is a subsolution of \eqref{14} lies below $(U_2(x),V_2(x,y))$ for all $h\in(-k,k)$. Anyway, by the uniqueness of the Cauchy problem associated to \eqref{14}, the solution of \eqref{14} with the translated subsolution as initial datum converges to the corresponding translation of $(U_2,V_2)$ and, by comparison, we have that $(U_2,V_2)$ is smaller than small translations in the $x$ direction of itself, which entails that the partial derivatives of $(U_2,V_2)$ with respect to $x$ are $0$.
\end{proof}
\end{proposition}

We are now able to give the

\vspace{0.2cm}
\noindent\emph{Proof of Theorem \ref{thltb}.} If \eqref{15} holds, we obtain \eqref{16} from Propositions \ref{pr31} and \ref{pr32}, since Proposition \ref{prcp1} guarantees that $(u,v)\geq(0,0)$. On the other hand, if \eqref{17} and \eqref{18} hold, \eqref{19} follows from Propositions \ref{pr31}, \ref{pr32} and \ref{pr34}. \hfill \qedsymbol

\vspace{0.2cm}
Since we are interested in the speed of propagation towards positive steady states of \eqref{14}, we will assume \eqref{17} and \eqref{18} throughout the rest of the paper, for \eqref{19} to hold.

\setcounter{equation}{0}
\setcounter{figure}{0}
\section{Supersolutions in the moving fra\-me\-work}
\label{sec4}
In this section we construct positive supersolutions to \eqref{14} moving at appropriate speeds in the $x-$direction. This will be the key to find an upper bound for the asymptotic speed of propagation of Theorem \ref{thmain} (see Section \ref{sec6}).
Observe that solutions of the linearized problem
\begin{equation}
\label{41}
  \left\{ \begin{array}{lll}
  u_t(x,t)-Du_{xx}(x,t)=\nu v(x,L,t)-\m u(x,t) & \text{for } x\in\R, &\!\! t>0 \\
  v_t(x,y,t)-d\D v(x,y,t)=f'(0)v(x,y,t) & \text{for } (x,y)\in\O, &\!\! t>0  \\ 
  dv_y(x,L,t)=\m u(x,t)-\nu v(x,L,t) & \text{for } x\in\R, &\!\! t>0 \\
  v(x,0,t)=0 & \text{for } x\in\R, &\!\! t>0, \end{array}\right.
\end{equation}
provide us with supersolutions to \eqref{14}, thanks to the second condition of \eqref{12}.
We start looking for solutions of \eqref{41} of the form
\begin{equation}
\label{42}
(u(x,t),v(x,y,t))=e^{\a(x+ct)}(1,\g \sin (\b y))
\end{equation}
with positive $\a,\g, c$ and $\b\in(0,\frac{\pi}{L})$, in order for $u$ and $v$ to be positive in $\O$.
Notice that \eqref{42} is a solution of \eqref{41} if and only if
\begin{equation}
\label{43}
  \left\{ \begin{array}{l}
  -D\a^2+c\a=\n\g\sin(\b L)-\m \\
  -d\a^2+d\b^2+c\a=f'(0) \\
  \g=\frac{\m}{d\b \cos(\b L)+\n \sin (\b L)}. \end{array}\right.
\end{equation}
In order for $\g$ to be positive, $\b$ must lie in $(0,\ov{\b})$, where $\ov{\b}\in(\frac{\pi}{2L},\frac{\pi}{L})$ is the first positive value of $\b$ for which
\begin{equation}
\label{44}
d\b\cos(\b L)+\n \sin(\b L)=0. 
\end{equation}

By substituting the expression of $\g$ given by the third equation of \eqref{43} into the first one, we get
\begin{equation*}
-D\a^2+c\a+\frac{\m}{1+\frac{\n\tan(\b L)}{d\b}}=0,
\end{equation*}
whose solutions are
\begin{equation}
\label{eqdef}
\a_D^{\pm}(c,\b)=\frac{1}{2D}\left(c\pm \sqrt{c^2+\chi(\b)}\right),
\end{equation}
where we have set
\begin{equation*}
\chi(\b)=\frac{4\m D}{1+\frac{\n\tan(\b L)}{d\b}}=\frac{4\m dD\b\cos(\b L)}{d\b\cos(\b L)+\n \sin(\b L)}.
\end{equation*}
It is easy to see that $\chi(\b)$ is continuous, even, decreases for $\b\in(0,\ov{\b})$ and satisfies
\begin{equation*}
\lim_{\b\da 0} \chi(\b)=\frac{4\m D}{1+\frac{\n L}{d}}, \qquad \chi\left(\frac{\pi}{2L}\right)=0, \qquad \lim_{\b\ua\ov{\b}} \chi(\b)=-\infty,
\end{equation*}
where $\ov{\b}$ is the one defined in \eqref{44}. Hence, for every $c>0$, there exists a unique $\tilde{\b}(c)\in(\frac{\pi}{2L},\ov{\b})$, for which
\begin{equation}
\label{45}
c^2=-\chi(\tilde{\b}(c)).
\end{equation}
Moreover it satisfies
\begin{equation*}
\lim_{c\da 0} \tilde{\b}(c)=\frac{\pi}{2L}, \qquad \lim_{c\ua \infty} \tilde{\b}(c)=\ov{\b}.
\end{equation*}

\begin{figure}[ht]
\begin{center}
\includegraphics[scale=0.76]{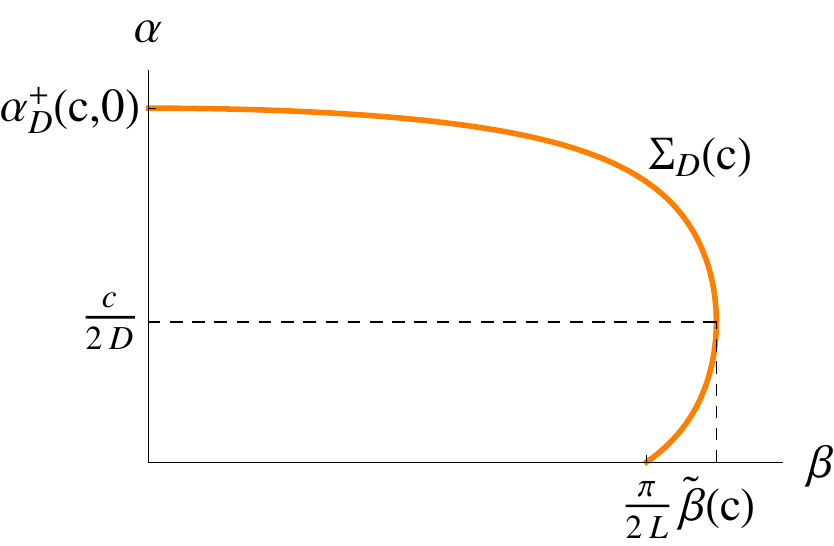}
\caption{The curve $\S_D(c)$ defined in \eqref{47}.} \label{fig21}
\end{center}
\end{figure}

As a consequence of these properties, for fixed $c>0$, $\a_D^{+}(c,\b)$ is a regular even function, which decreases in $(0,\tilde{\b}(c))$ and satisfies
\begin{gather}
\a_D^+(c,0)\!=\!\frac{1}{2D}\left(c\!+\!\sqrt{c^2+\frac{4\m D}{1+\frac{\n L}{d}}}\right), \quad \a_D^{+}\left(c,\frac{\pi}{2L}\right)\!=\!\frac{c}{D}, \quad \a_D^{+}(c,\tilde{\b}(c))\!=\!\frac{c}{2D}, \label{eqstar} \\ 
\p_\b\a_D^+(c,0)=0, \qquad \lim_{\b\ua\tilde{\b}(c)} \p_\b\a_D^+(c,\b)=-\infty. \notag
\end{gather}
In addition it is easy to verify that
\begin{equation}
\label{46}
\p_c\a_D^+(c,\b)>0. 
\end{equation}
Since $\a_D^-$ is symmetric to $\a_D^+$ with respect to the line $\a=\frac{c}{2D}$, analogous properties can be established for $\a_D^-$. In particular observe that
\begin{equation*}
\a_D^{-}\left(c,\frac{\pi}{2L}\right)=0
\end{equation*}
and that, for fixed $\b\in(\frac{\pi}{2L},\ov{\b})$, we have
\begin{align*}
\lim_{c\ua\infty} \a_D^{-}(c,\b)=0,
\end{align*}
$\a_D^{-}$ decreasing monotonically in $c$.
As a consequence, the bounded region of the first quadrant in the $(\b,\a)-$plane delimited by the curve
\begin{equation}
\label{47}
\S_D(c):=\left\{(\b,\a_D^-(\b)):\b\in\left[\frac{\pi}{2L},\tilde{\b}(c)\right]\right\}\cup\left\{(\b,\a_D^+(\b)):\b\in\left(0,\tilde{\b}(c)\right]\right\},
\end{equation}
which has been represented in Figure \ref{fig21}, invades monotonically, as $c\ua\infty$, the half-strip $\{(\b,\a): \b\in(0,\ov{\b}), \a>0\}$.

As far as the monotonicity of the curves with respect to $D$, the other parameters being fixed, it follows from the definition of $\a_D^{\pm}$ that
\begin{equation}
\label{48}
\begin{split}
D&\mapsto \a_D^+(c,\b) \; \text{ is decreasing for every $\b\in[0,\tilde{\b}(c)]$} \\
D&\mapsto \a_D^-(c,\b) \; \text{ is increasing for every $\b\in\left(\frac{\pi}{2L},\tilde{\b}(c)\right]$.}
\end{split}
\end{equation}

\begin{figure}[ht]
\begin{center}
\begin{tabular}{ccc}
\includegraphics[scale=0.45]{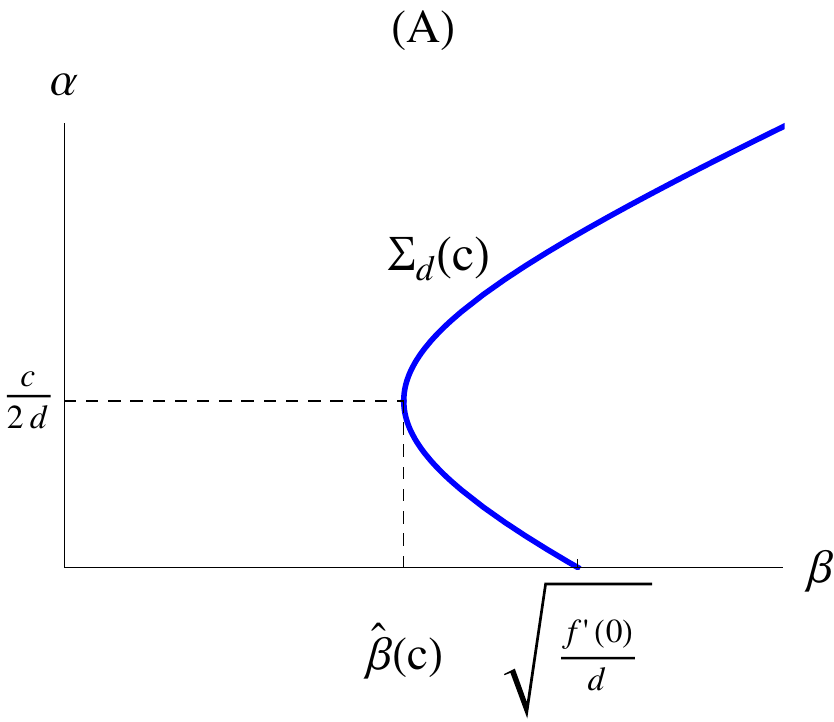} & \hspace{0.1cm}
\includegraphics[scale=0.46]{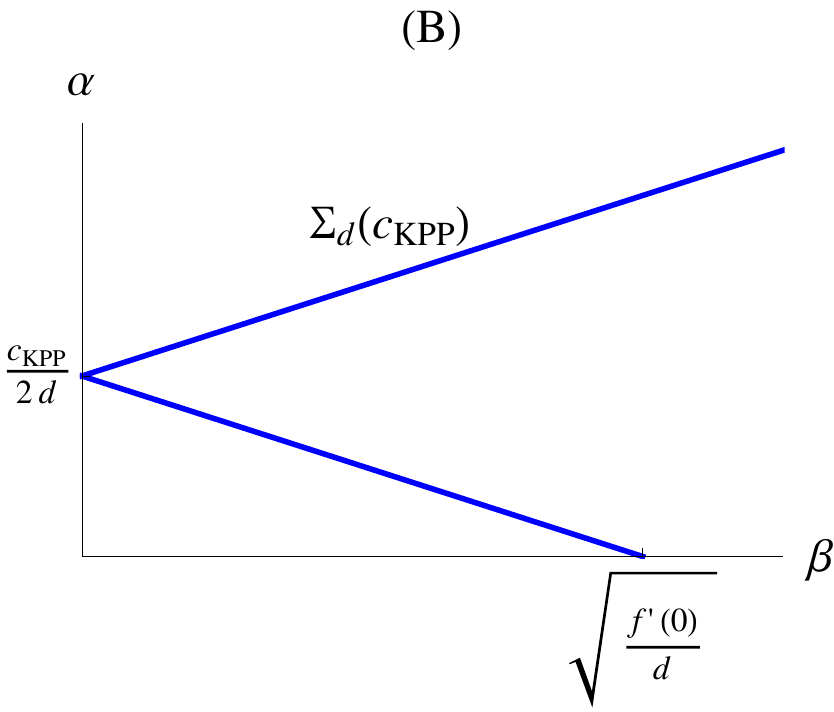} & \hspace{0.1cm}
\includegraphics[scale=0.49]{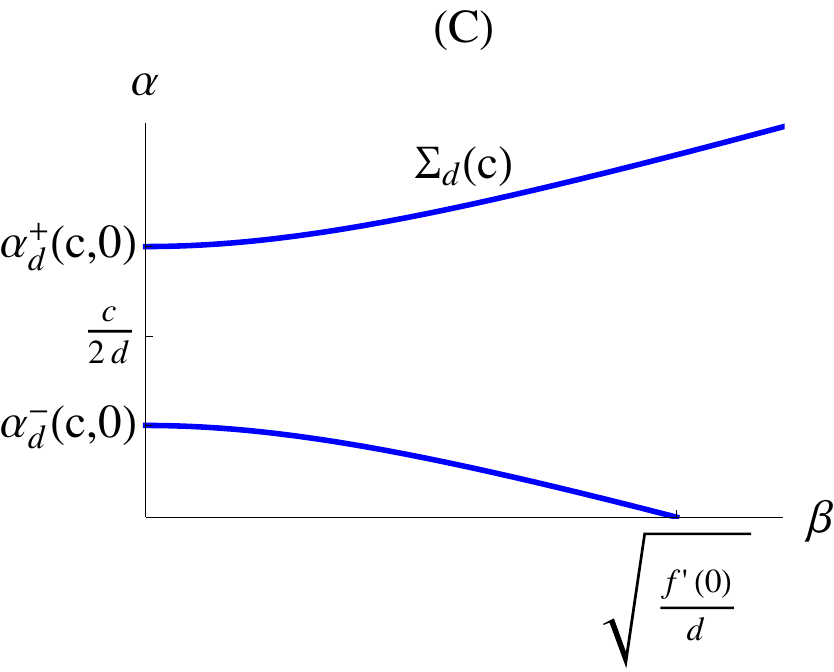} 
\end{tabular}
\caption{The hyperbola defined by the second equation of \eqref{43}, for $0<c<\cKPP$ (A), $c=\cKPP$ (B) and $c>\cKPP$ (C).} \label{fig22}
\end{center}
\end{figure}

On the other hand (see Figure \ref{fig22}), the second equation of \eqref{43} represents an hyperbola whose branches are given by
\begin{equation}
\label{49}
\a_d^{\pm}(c,\b)=\frac{1}{2d}\left(c\pm\sqrt{c^2-\eta(\b)}\right),
\end{equation}
where we have set,
\begin{equation*}
\eta(\b):=4d(f'(0)-d\b^2)=\cKPP^2-4d^2\b^2,
\end{equation*}
$c_{\KPP}$ being the one defined in \eqref{13}. The function $\eta$ is decreasing and satisfies
\begin{equation*}
\eta(0)=\cKPP^2, \qquad \eta\left(\sqrt{\frac{f'(0)}{d}}\right)=0, \qquad \lim_{\b\ua\infty}\eta(\b)=-\infty.
\end{equation*} 
Therefore, the functions $\a_d^{\pm}(c,\b)$ are defined for every $\b\in(0,+\infty)$ if $c\geq \cKPP$, while, if $c\in(0,\cKPP)$ there exists $\hat{\b}(c)>0$ such that their domain (within the positive part of the real line) is $[\hat{\b}(c),+\infty)$, where $\hat{\b}(c)$ satisfies
\begin{equation}
\label{410}
c^2=\eta(\hat{\b}(c)).
\end{equation}
As a consequence,
\begin{equation*}
\lim_{c\da 0} \hat{\b}(c)=\sqrt{\frac{f'(0)}{d}}, \qquad \lim_{c\ua \cKPP} \hat{\b}(c)=0.
\end{equation*}
It can be easily seen, as $\a_d^-$ is the reflection of $\a_d^+$ about the line $\a=\frac{c}{2d}$, that
\begin{equation*}
\a_d^-\left(c,\sqrt{\frac{f'(0)}{d}}\right)=0, \qquad \lim_{\b\ua\infty}\a_d^+(c,\b)=+\infty,
\end{equation*}
and, in the proper domain of definition,
\begin{equation}
\label{411}
\p_\b\a_d^+(c,\b)>0, \qquad \p_c\a_d^+(c,\b)>0, \qquad \p_c\a_d^-(c,\b)<0.
\end{equation}
Moreover, if $0<c<\cKPP$, we have that
\begin{equation*}
\a_d^+(c,\hat{\b}(c))=\frac{c}{2d} \qquad \lim_{\b\da\hat{\b}(c)}\p_\b\a_d^{+}(c,\b)=+\infty,
\end{equation*}
while, for $c=\cKPP$, the hyperbolas degenerate into the straight lines with equations
\begin{equation}
\label{ipdeg}
\a_d^{\pm}(\cKPP,\b)=\pm\b+\frac{\cKPP}{2d}.
\end{equation}
Finally, for $c>\cKPP$, we have that
\begin{equation*}
\a_d^{\pm}(c,0)=\frac{c\pm\sqrt{c^2-\cKPP^2}}{2d}, \qquad \p_\b\a_d^{\pm}(c,0)=0
\end{equation*}
and, for fixed $\b$,
\begin{equation*}
\lim_{c\ua\infty} \a_d^{-}(c,\b)=0.
\end{equation*}
As a consequence of all the aforementioned properties, if we set $\hat{\b}(c)=0$ for $c\geq \cKPP$, which is consistent with the previous notation, and define
\begin{equation*}
\S_d(c):=\left\{(\b,\a_d^-(\b)):\b\in\left[\hat{\b}(c),\sqrt{\frac{f'(0)}{d}}\right]\right\}\cup\left\{(\b,\a_d^+(\b)):\b\geq\hat{\b}(c)\right\},
\end{equation*}
we have that the region of the first quadrant in the $(\b,\a)-$plane bounded by $\S_d(c)$ and containing the point $\left(\sqrt{\frac{f'(0)}{d}},\frac{c}{2d}\right)$
invades monotonically, as $c\ua\infty$, the first quadrant in the $(\b,\a)-$plane. All these features have been represented in Figure \ref{fig22}.

Another candidate to construct supersolutions to \eqref{14} for $c>\cKPP$ are the functions
\begin{equation}
\label{412}
(u(x,t),v(x,y,t))=e^{\a(x+ct)}(1,\phi(y))
\end{equation}
with
\begin{equation*}
\phi(y):=\g(e^{\b y}-e^{-\b y}),
\end{equation*}
where $\a,\b,\g$ are positive constants. By plugging \eqref{412} into \eqref{41} we are driven to
\begin{equation}
\label{413}
  \left\{ \begin{array}{l}
  -D\a^2+c\a=\n\phi(L)-\m \\
  -d\a^2-d\b^2+c\a=f'(0) \\
  \g=\frac{\m}{d\b(e^{\b L}+e^{-\b L})+\n (e^{\b L}-e^{-\b L})}>0. \end{array}\right.
\end{equation}
Using the expression of $\phi(L)$ given from the third equation of \eqref{413}, the first one reduces to
\begin{equation*}
-D\a^2+c\a+\frac{\m d\b(e^{\b L}+e^{-\b L})}{d\b(e^{\b L}+e^{-\b L})+\n (e^{\b L}-e^{-\b L})}=0,
\end{equation*}
whose solutions are
\begin{equation}
\label{tipo2}
\tilde{\a}_D^{\pm}(c,\b)=\frac{1}{2D}\left(c\pm \sqrt{c^2+\tilde{\chi}(\b)}\right),
\end{equation}
where we have set
\begin{equation*}
\tilde{\chi}(\b)=\frac{4\m D}{1+\frac{\n\tanh(\b L)}{d\b}}.
\end{equation*}
This function is positive and therefore $\tilde{\a}_D^{+}(c,\b)$ is defined for every $\b\in\R$. It can be easily seen that it is even and satisfies the following monotonicity conditions
\begin{equation}
\label{417}
\p_{\b}\tilde{\a}_D^+(c,\b)>0, \quad \p_{c}\tilde{\a}_D^+(c,\b)>0, \quad \p_{D}\tilde{\a}_D^+(c,\b)<0
\end{equation}
for every positive $c,\b$. Moreover we have $\tilde{\a}_D^+(c,0)=\a_D^+(c,0)$,
for every $c>0$. Naturally, similar properties hold for $\tilde{\a}_D^-$, taking into account that it is symmetric to $\tilde{\a}_D^+$ with respect to $\a=\frac{c}{2D}$. Finally, observe that, for every $\b\geq 0$,
\begin{equation*}
\lim_{c\ua+\infty}\tilde{\a}_D^+(c,\b)=+\infty
\end{equation*}
and, therefore, the regions of the first quadrant in the $(\b,\a)-$plane delimited by
\begin{equation}
\label{418}
\tilde{\S}_D(c):=\left\{(\b,\tilde{a}_D^-(c,\b)):\b>0\right\}\cup\left\{(\b,\tilde{a}_D^+(c,\b)):\b>0\right\}.
\end{equation}
and containing $\left(\frac{c}{2D},1\right)$
invade monotonically the first quadrant of the plane $(\b,\a)$.

\begin{figure}[ht]
\begin{center}
\includegraphics[scale=1.0]{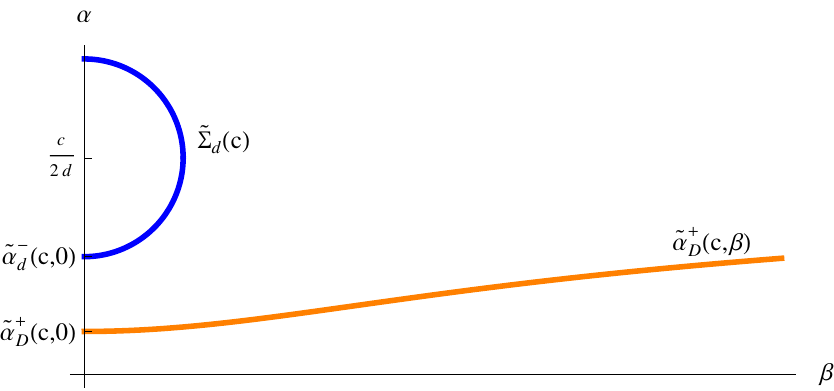}
\caption{The curves $\tilde{\a}_D^+$ and $\tilde{\a}_d^{\pm}$.} \label{fig24}
\end{center}
\end{figure}

On the other hand, the second equation of \eqref{413} describes, for $c>\cKPP$, as we are assuming, a circle in the $(\b,\a)-$plane, with center at $\left(0,\frac{c}{2d}\right)$ and radius
\begin{equation*}
r(c)=\frac{\sqrt{c^2-\cKPP^2}}{2d}
\end{equation*}
(see Figure \ref{fig24}). Precisely, the part of the graph of the circle which lies in the first quadrant is given by
\begin{equation*}
\tilde{\S}_d(c):=\left\{(\b,\tilde{\a}_d^{\pm}(c,\b)):\b\in[0,r(c)]\right\},
\end{equation*}
where
\begin{equation}
\label{419}
\tilde{\a}_d^{\pm}(c,\b)=\frac{1}{2d}\left(c\pm\sqrt{c^2-\cKPP^2-4d^2\b^2}\right).
\end{equation}
The function $\tilde{\a}_d^-$ satisfies $\p_{\b}\tilde{\a}_d^-(c,\b)>0$, $\p_{c}\tilde{\a}_d^-(c,\b)<0$, for $c>\cKPP$ and $\b\in(0,r(c))$, while its value 
at $\b=0$ satisfies $\tilde{\a}_d^-(c,0)=\a_d^-(c,0)$.
Moreover, it can be easily seen that
\begin{equation*}
\tilde{\S}_d(\cKPP)=\left\{\left(0,\frac{\cKPP}{2d}\right)\right\}
\end{equation*}
and
\begin{equation*}
\lim_{c\ua+\infty}\tilde{a}_d^-(c,\b)=0.
\end{equation*}
As a consequence of these properties, the half-disks delimited by $\tilde{\S}_d(c)$ and contained in the first quadrant of the $(\b,\a)-$plane invade it monotonically as $c$ increases.
With these ingredients we are able to give the following result

\begin{proposition}
\label{pr41}
There exists $c^*:=c^*_L(D,d,\mu,\nu)$ such that, for every $c>c^*$, Problem \eqref{14} admits supersolutions either of the form \eqref{42} or \eqref{412}, with positive $\a,\b,\g$.
Moreover the function $D\mapsto c^*(D)$ is increasing.

\begin{proof}
By the previous discussion, we have that, in order to find a solution of \eqref{41} and therefore a supersolution to \eqref{14}, it is sufficient to find an intersection between either $\S_D(c)$ and $\S_d(c)$ or $\tilde{\S}_D(c)$ and $\tilde{\S}_d(c)$ lying in the interior of the first quadrant of the $(\b,\a)-$plane. Due to the monotonicity properties of these curves with respect to $c$ shown above, $c^*$ will be the smallest value of $c$ for which such an intersection exists for all $c>c^*$.

Let us start by examining the case $D$ small (relatively to $d$), in which we consider the curves $\S_D(c)$ and $\S_d(c)$. Recalling that we are assuming \eqref{17}, it is clear from the above discussion that, for $c\sim 0$, they are disjoint, since so are their domain of definition. On the contrary, for sufficiently large $c$ such an intersection exists, since
\begin{equation*}
\a_d^-(c,0)<\a_D^+(c,0), \quad \a_d^-\left(c,\sqrt{\frac{f'(0)}{d}}\right)=0, \quad \a_D^-\left(c,\frac{\pi}{2L}\right)=0.
\end{equation*}

\begin{figure}[ht]
\begin{center}
\begin{tabular}{cc}
\includegraphics[scale=0.64]{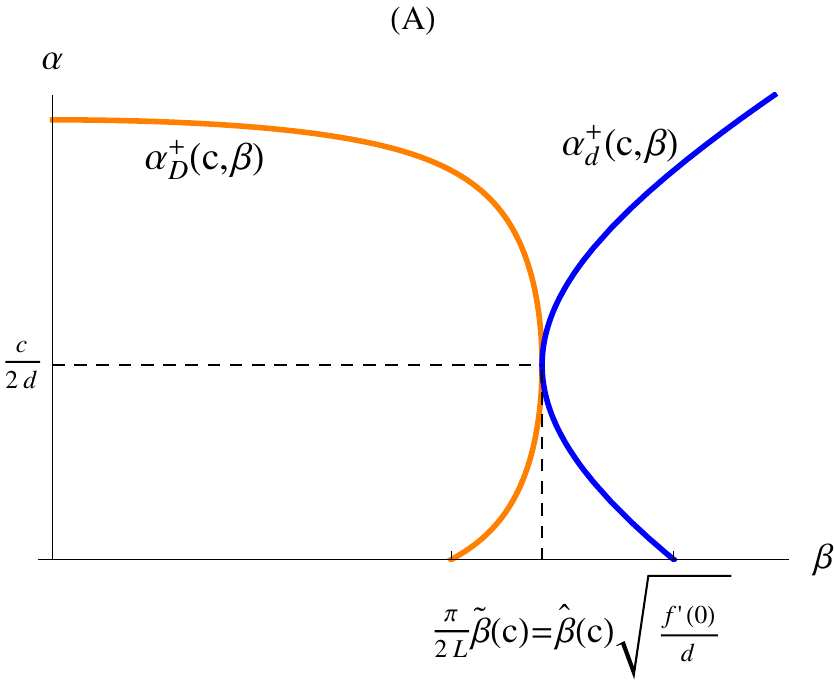} & \hspace{0.1cm}
\includegraphics[scale=0.64]{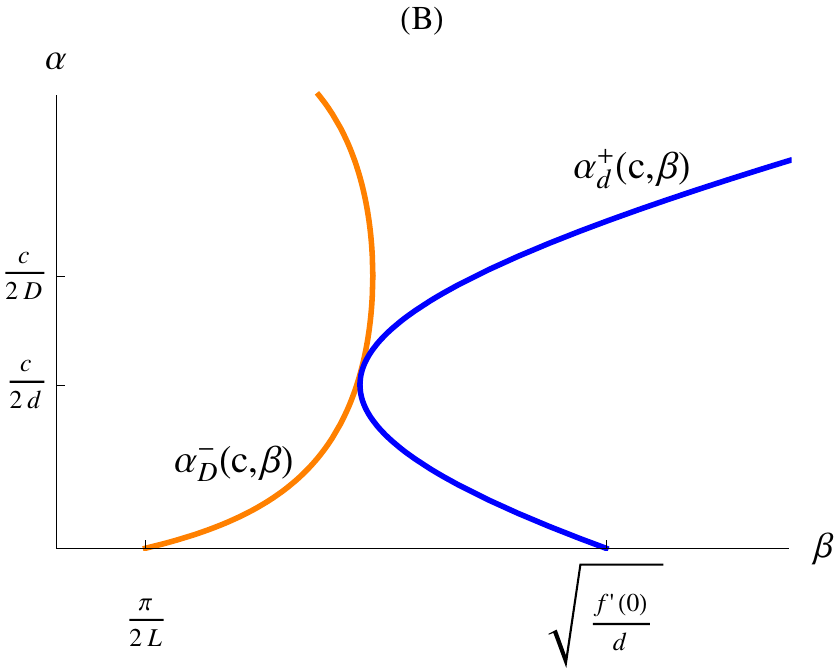} \\
\includegraphics[scale=0.64]{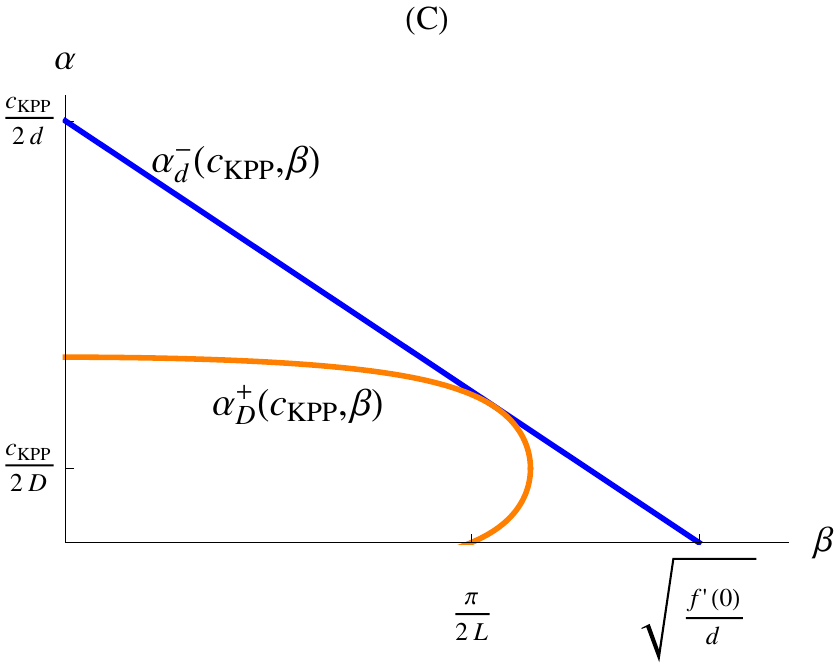} & \hspace{0.1cm}
\includegraphics[scale=0.64]{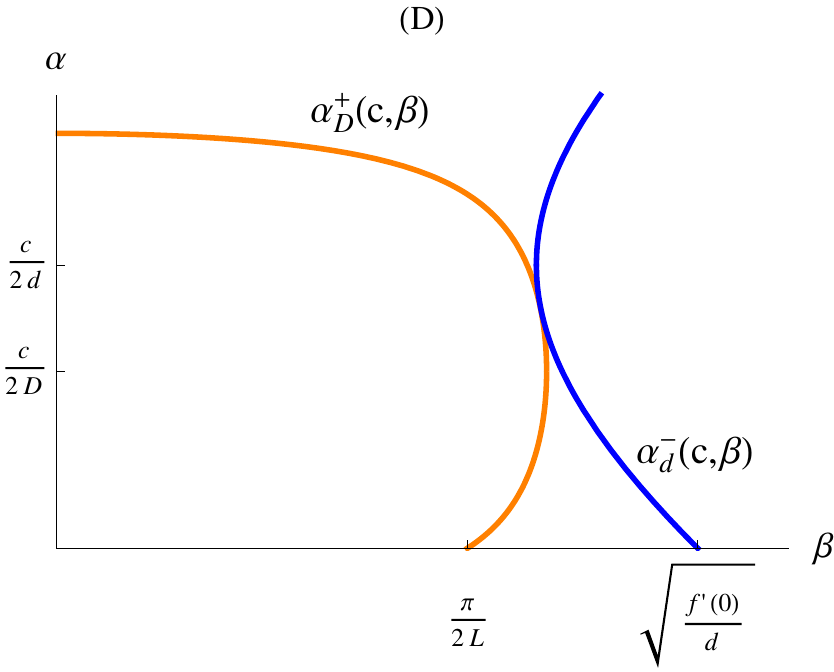} \\
\end{tabular}
\caption{The tangency between $\S_D(c)$ and $\S_d(c)$, for $D\leq D_{\KPP}$.} \label{fig23}
\end{center}
\end{figure}

Actually, the first value of $c$, denoted by $\un{c}=\un{c}(D)$, for which an intersection can exist is the one for which the domains of $\a_D^{\pm}$ and $\a_d^{\pm}$ intersect in one point, i.e. the one for which
\begin{equation*}
  \left\{ \begin{array}{l}
  c^2=-\chi(\b) \\
  c^2=\eta(\b) \end{array}\right.
\end{equation*}
admits a solution. Indeed, as shown in Figure \ref{fig23}(A), the curves $\S_D(\un{c})$ and $\S_d(\un{c})$ are tangent if $D=d$, their common tangent being vertical. Therefore in this case $c^*=\un{c}$. If $D<d$,  there exists $c^*$ satisfying $\un{c}<c^*<\cKPP$ and such that $\a_D^-$ and $\a_d^+$ are tangent, as described in Figure \ref{fig23}(B). If $D>d$ the situation is more complex, because we have to take into account the change in the nature of $\a_d^-(c,\b)$ as $c$ crosses $\cKPP$. From \eqref{45}, it follows that
\begin{equation*}
\lim_{D\ua+\infty}\tilde{\b}(c_{\KPP},D)=\frac{\pi}{2L}
\end{equation*}
and, together with the first and third relations of \eqref{eqstar}, this implies that there exists a value of $D$, denoted by $D_{\KPP}$, for which, for $c=c_{\KPP}$, $\a_{D_{\KPP}}^+$ and the straight line $\a_d^-$ are tangent (see Figure \ref{fig23}(C)). Observe that, for $D=2d$, we have 
\begin{equation*}
\a_{D}^+(c_{\KPP},0)>\a_{D}^+\left(c_{\KPP},\frac{\pi}{2L}\right)=\frac{c_{\KPP}}{D}=\frac{c_{\KPP}}{2d}=\a_{d}^-(c_{\KPP},0),
\end{equation*}
which implies that $D_{\KPP}>2d$. Thanks to \eqref{48}, for $d<D<D_{\KPP}$ and $c=c_{\KPP}$, $\a_D^+$ and $\a_d^-$ will be secant, therefore in this case the tangency will occur between $\a_D^+$ and $\a_d^-$ for $\un{c}<c^*<\cKPP$, as represented in Figure \ref{fig23}(D).

\begin{figure}[ht]
\begin{center}
\begin{tabular}{ccc}
\includegraphics[scale=0.48]{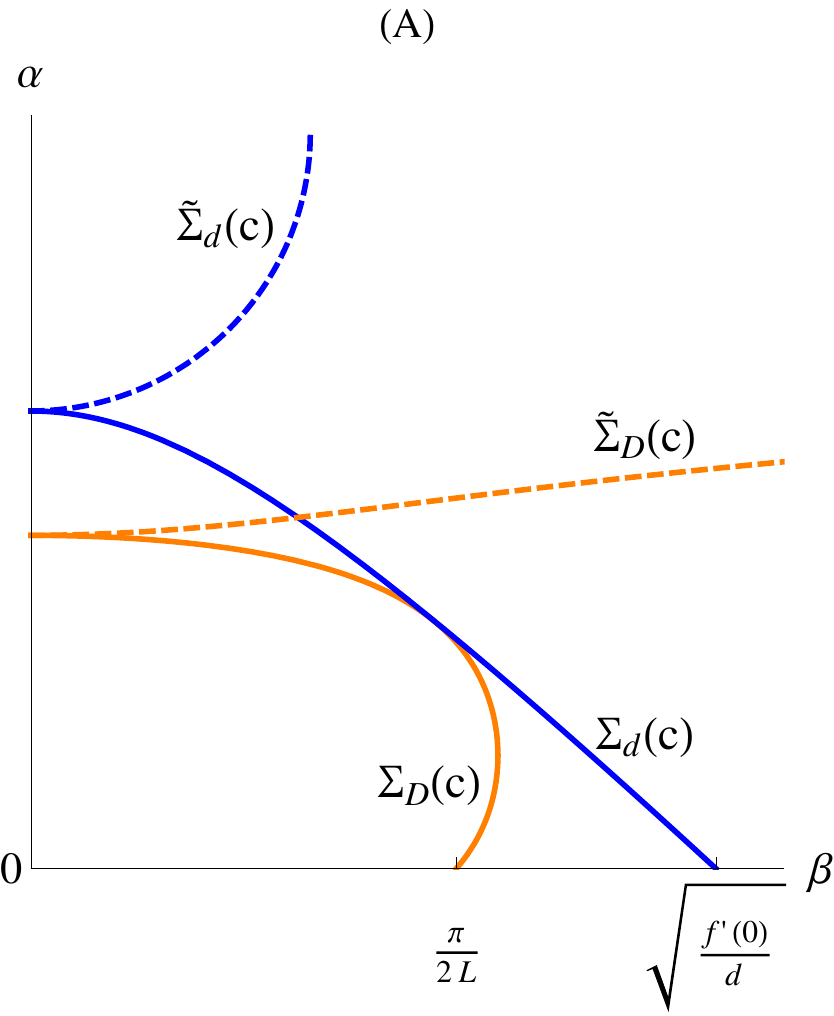} & \hspace{0.1cm}
\includegraphics[scale=0.48]{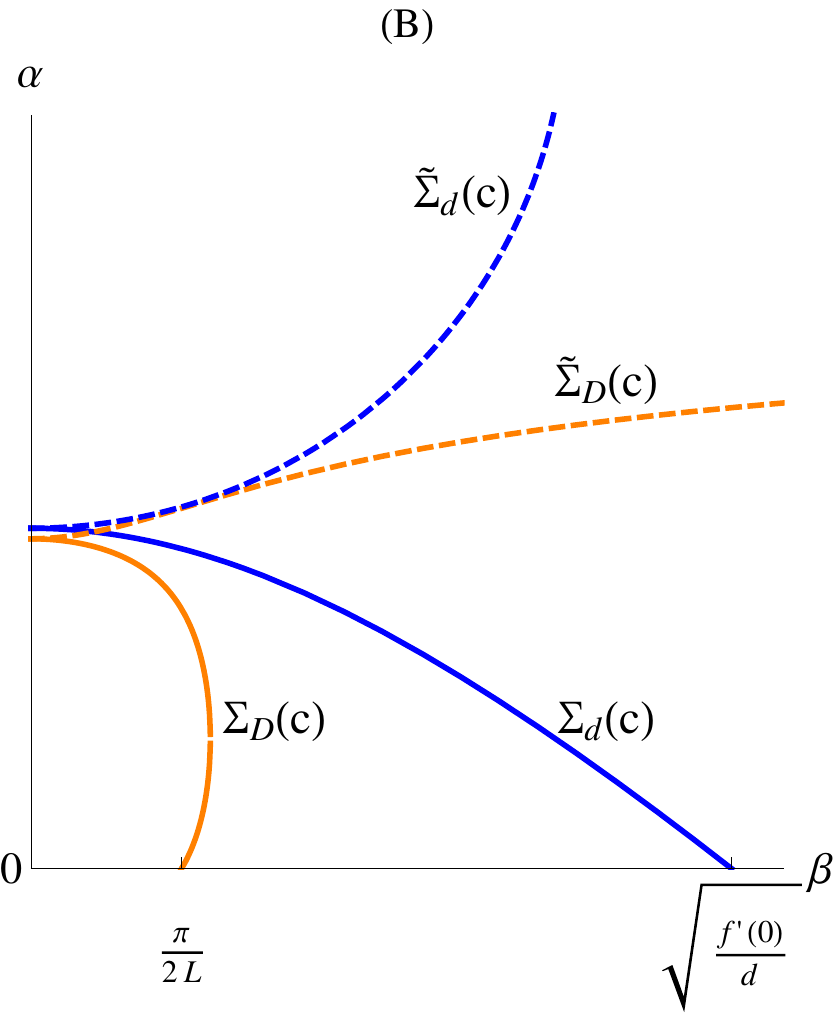} & \hspace{0.1cm}
\includegraphics[scale=0.48]{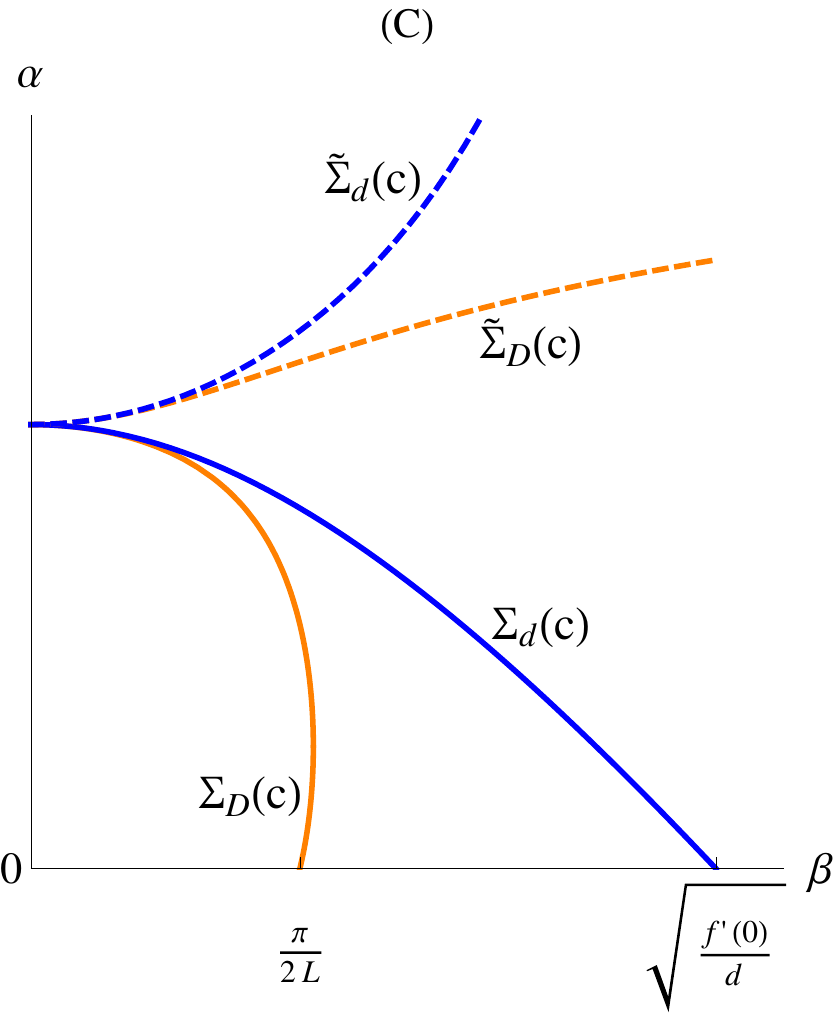}
\end{tabular}
\caption{Different configurations of the tangency between $\S_D(c)$ and $\S_d(c)$ (continuous lines) or  $\tilde{\S}_D(c)$ and $\tilde{\S}_d(c)$ (dashed lines), for $c>c_{\KPP}$.} \label{fig26}
\end{center}
\end{figure}

Finally, when $D>D_{\KPP}$ we consider at the same time the pairs of curves $\S_D(c)$, $\S_d(c)$ and $\tilde{\S}_D(c)$, $\tilde{\S}_d(c)$. We define $c^*$ to be the smallest value of $c$ for which either the former two curves or the latter two are tangent in a positive $\b$ (see Figure \ref{fig26}(A) and (B) respectively). The last case to be examined is the one in which all the four curves touch for the first time, being tangent, at $\b=0$ for $c=\cint$, where $\cint$ is
\begin{equation}
\label{420}
\cint^2=\cint^2(D)=\frac{\left((d+\nu L)D\cKPP^2+ 4\mu d^3\right)^2}{4d(D-d)(d+\nu L)\left((d+\nu L)\cKPP^2 + 4\mu d^2\right)}
\end{equation}
(see figure \ref{fig26}(C)). In this case, thanks to \eqref{46} and \eqref{411}, we have that $\S_D(c)$ and $\S_d(c)$ intersect for every $c>\cint$ and a certain $\b>0$. Therefore we set $c^*=\cint$.

The monotonicity of the function $D\mapsto c^*(D)$ follows from the monotonicity of the curves $\a_D^{\pm}$ and $\tilde{\a}_D^+$, given by \eqref{48} and the last relation of \eqref{417} respectively.
\end{proof}
\end{proposition}

\begin{remark}
\label{remark42}
\begin{itemize}
\item[(i)] The proof of Proposition \ref{pr41} shows indeed that it is possible to construct a supersolution to \eqref{14} of type \eqref{42} or \eqref{412} not only for every $c>c^*$, but also for $c=c^*$, except for the case represented in Figure \ref{fig26}(C). Actually it is possible to construct a supersolution also in this case, when $c^*=\cint$, by taking
\begin{equation}
\label{421}
(u(x,t),v(x,y,t))=e^{\a(x+ct)}\left(1,\frac{\mu y}{d+\nu L}\right).
\end{equation}
This can be heuristically seen by taking the limit of \eqref{42} or \eqref{412} as $\b\da 0$ (notice that $\g=\g(\b)\to 0$), while a formal proof consists in plugging \eqref{421} into \eqref{41} and observing that the resulting algebraic system in $\a,c$ has a solution for $c=\cint$.
\item[(ii)] From the last part of proof of Proposition \ref{pr41}, it arises the natural question of characterizing which of the curves first touch, either $\S_D(c)$ and $\S_d(c)$ or $\tilde\S_D(c)$ and $\tilde\S_d(c)$. This analysis can be performed following the same ideas of the proof of Theorem \ref{thmain}(iii) (see Section \ref{sec8}), i.e. by studying, for $c=\cint$, the second derivatives of the curves at $\b=0$, and adapting a fine result given in \cite[Proposition 4.1]{RTV}, which is based on the comparison principles and characterizes the total number of possible intersections between the curves. Since such analysis is quite technical in general (indeed some of the arguments provided in Section \ref{sec8} do not adapt to the general situation) and not strictly necessary for the results of this work, we send the reader to \cite{RTV} for the details.
\end{itemize}
\end{remark}

\setcounter{equation}{0}
\setcounter{figure}{0}
\section{Generalized subsolutions with compact support}
\label{sec5}
In this section we construct stationary compactly supported generalized (in the sense of Proposition \ref{prcp2}) subsolutions in a framework moving in the $x-$direction at slightly smaller speeds than $c^*$, the one of Proposition \ref{pr41}. Provided that Proposition \ref{prcp2} can be applied, this will provide a lower bound for the asymptotic speed of propagation and will be the second and last ingredient for the proof of Theorem \ref{thmain} (see Section \ref{sec6}). The result is the following
\begin{proposition}
\label{pr51}
Let $c^*$ be the one constructed in Proposition \ref{pr41}. Then, for every $c<c^*$, $c\sim c^*$, and $\d>0$, $\d\sim 0$, there exists $\e_0>0$ such that, for every $0<\e<\e_0$, $\e(\un{U}(x),\un{V}(x,y))$
is a compactly supported generalized (in the sense of Proposition \ref{prcp2}) subsolution of
\begin{equation}
\label{51}
  \left\{ \begin{array}{lll}
  u_t-Du_{xx}+c u_x=\nu v(x,L,t)-\m u & \text{for } x\in\R, & t>0 \\
  v_t-d\D v+c v_x=f(v) & \text{for } (x,y)\in\O, & t>0  \\ 
  dv_y(x,L,t)=\m u-\nu v(x,L,t) & \text{for } x\in\R, & t>0 \\
  v(x,0,t)=0 & \text{for } x\in\R, & t>0, \end{array}\right.
\end{equation}
satisfying \eqref{21}, where $(\un{U}(x),\un{V}(x,y))$ is a truncation of a solution of
\begin{equation}
\label{52}
  \left\{ \begin{array}{ll}
  -D U_{xx}(x)+c U_x(x)=\nu V(x,L)-\m U(x) & \text{for } x\in\R, \\
  -d\D V(x,y)+c V_x(x,y)=(f'(0)-\d)V(x,y) & \text{for } (x,y)\in\O,   \\ 
  d V_y(x,L)=\m U(x)-\nu V(x,L) & \text{for } x\in\R,  \\
  V(x,0)=0 & \text{for } x\in\R.  \end{array}\right.
\end{equation}

\begin{proof}
Once  the solution $(\un{U},\un{V})$ of \eqref{52} will be constructed, the existence of $\e_0$ such that $\e(\un{U}(x),\un{V}(x,y))$ is a subsolution to \eqref{51} for every $0<\e<\e_0$ follows immediately from \eqref{12}. From the construction of $c^*$ carried out in the proof of Proposition \ref{pr41}, we know that, for $c<c^*$, no real solution of \eqref{52} of type \eqref{42} or \eqref{412} exists. However we will show that \eqref{52} exhibits, for $c<c^*$, $c\sim c^*$, complex solutions, which will be the starting point for the construction of $(\un U,\un V)$. Actually we will consider the case $\d=0$ and the existence of $(\un{U},\un{V})$ for $\d\sim 0$ will follow by a perturbation argument which consists in repeating the same construction of the proof of Proposition \ref{pr41} considering the curves with $f'(0)$ replaced by $f'(0)-\d$ (observe that the dependence of the curves with respect to $f'(0)$ is continuous).

We will give the details only in the case in which $c^*$ was constructed in Proposition \ref{pr41} as the one for which $\a_D^+(c^*,\b)$ and $\a_d^-(c^*,\b)$ were tangent in a point $\b=\b^*>0$. The other cases related to supersolution like \eqref{42} are analogous; the case of supersolutions like \eqref{412} was treated in \cite{BRR1} and the one related to supersolutions like \eqref{421} will follow from the case that we are going to present here, by passing to the limit as $\b^*\to 0$, like in Remark \ref{remark42}(i). 

Let us consider, for $(c,\b)$ in a neighborhood of $(c^*,\b^*)$, the function
\begin{equation}
\label{53}
g(c,\b)=\a_d^-(c,\b)-\a_D^+(c,\b).
\end{equation}
Our goal is to find, for $c<c^*$, $c\sim c^*$, a root $\b\in\C\setminus \R$ of \eqref{53}. In this way we will obtain a solution $(\a,\b,\g)=(\a_1+i\a_2,\b_1+i\b_2,\g_1+i\g_2)\in(\C\setminus \R)^3$ of \eqref{43}. It is easily seen that $(\ov{\a},\ov{\b},\ov{\g})$ also solves \eqref{43} and, therefore, by taking the real part in \eqref{42}, we can set
\begin{gather}
\label{54}
\un{U}(x)=\max\{e^{\a_1x}\cos(\a_2x), 0\}, \\
\begin{split}
\label{55}
\un{V}(x,y)=&\max\{e^{\a_1x}(\sin(\b_1 y)\cosh(\b_2 y)(\g_1\cos(\a_2x)-\g_2\sin(\a_2x))+ \\
&- \cos(\b_1 y)\sinh(\b_2 y)(\g_1\sin(\a_2x)+\g_2\cos(\a_2x))), 0\},
\end{split}
\end{gather}
where
\begin{equation}
\label{56}
\g_1=\frac{\mu\t_1}{\t_1^2+\t_2^2}, \qquad
\g_2=\frac{\mu\t_2}{\t_1^2+\t_2^2},
\end{equation}
being
\begin{equation}
\label{57}
\begin{split}
\t_1& = ( d\b_1\cos( \b_1 L ) + \nu\sin( \b_1 L ) )\cosh( \b_2 L )+d\b_2\sin( \b_1 L )\sinh( \b_2 L ) \\
\t_2& = ( d\b_1\sin( \b_1 L ) - \nu\cos( \b_1 L ) )\sinh( \b_2 L ) - d\b_2\cos( \b_1 L )\cosh( \b_2 L ).
\end{split}
\end{equation}
After the change of variables
\begin{equation*}
\xi=\b-\b^*, \qquad \tau=c-c^*
\end{equation*}
the search for zeros of \eqref{53} is equivalent to the search for zeros of the function
\begin{equation}
\label{58}
h(\xi,\tau):=g(c^*+\tau,\b^*+\xi)
\end{equation}
with $(\xi,\tau)$ in a neighborhood of $(0,0)$. Since $(c^*,\b^*)$ is the first contact point between $\a_D^+$ and $\a_d^-$, we have that there exists $n\in\N\setminus\{0\}$ such that
\begin{equation*}
h(0,0)=\dots=\p^{2n-1}_\xi h(0,0)=0 \quad \text{ and } a_{2n}:=\frac{\p^{2n}_\xi h(0,0)}{(2n)!}>0,
\end{equation*}
while from \eqref{46} and \eqref{411} it follows that
\begin{equation*}
a_1:=\p_\tau h(0,0)<0.
\end{equation*}
By considering the Taylor series of \eqref{58} in a neighborhood of $(0,0)$ we have that $h(\xi,\tau)=0$ is equivalent to
\begin{equation}
\label{59}
a_1\tau+a_{2n}\xi^{2n}=p(\xi,\tau)\tau+o(\xi^{2n+1})
\end{equation}  
where $p(\xi,\tau)$ is a polynomial which is either identically $0$ or of degree at least $1$. Thanks to the signs of the coefficients determined above, we know that the left hand side of \eqref{59}
\begin{equation*}
h_1(z):=a_{2n}z^{2n}+a_1\tau
\end{equation*}
has, for $\tau<0$, $2n$ complex roots
\begin{equation*}
z_j:=z_j(\tau)=\left(\frac{a_1 \tau}{a_{2n}}\right)^{\frac{1}{2n}}e^{i\frac{(2j-1)\pi}{2n}} \quad j=1,\dots,2n.
\end{equation*}
Consider now the ball
\begin{equation*}
B:=B_r(z_1)\subset\C \qquad  \text{ with } r=\s|\tau|^{\frac{1}{2n}}, \quad \s\sim 0
\end{equation*}
From geometrical considerations we have that, on $\p B$,
\begin{equation*}
|h_1(z)|=\prod_{j=1}^{2n}|z-z_j|\geq r\prod_{j=2}^{2n}C_j|\tau|^{\frac{1}{2n}}>C|\tau|
\end{equation*}
while, the right hand side $\varphi$ of \eqref{59}, considered as a function of $\xi$, satisfies, on $\p B$,
\begin{equation*}
|\varphi(z)|\leq|p(z,\tau)||\tau|+o(|z|^{2n+1})<\tilde{C}|\tau|^{1+\frac{1}{2n}}.
\end{equation*}
Therefore, by choosing $\tau$ negative and sufficiently small, we can make $|\varphi|<|h_1|$ on $\p B$ and Rouch\'e's theorem can be applied, guaranteeing the existence of complex roots of \eqref{58} and therefore of \eqref{53} for $c<c^*$, $c\sim c^*$. This same analysis also shows that $\b=\b_1+i\b_2=\b_1(c)+i\b_2(c)$ satisfies
\begin{equation}
\label{510}
\b_1(c)\to \b^*, \quad \b_2(c)\to 0 \qquad \text{ as } c\ua c^*.
\end{equation}
As a consequence, from \eqref{56} and  \eqref{57} we have that
\begin{equation*}
\cosh(\b_2 L)\to 1, \quad \sinh(\b_2 L)\to 0, \quad \g_1\to\frac{\mu}{d\b^*\cos(\b^*L)+\nu\sin(\b^*L)}>0, \quad \g_2\to 0,
\end{equation*}
as $c\ua c^*$, since $0<\b^*<\ov{\b}<\pi/L$. Moreover, when $n>1$, the second equations in systems \eqref{43} and \eqref{413} ensure that $\a_2\neq 0$, since both $\b_1$ and $\b_2$ are positive for $c\sim c^*$. This follows from \eqref{510} by continuity in the case $\b^*\neq 0$ and, when $\b^*=0$, it can be shown directly by analyzing the equation for $\a$ which is obtained by plugging \eqref{421} into \eqref{41}. On the other hand, when $n=1$, it can be proved as in Section \ref{sec8} (see \eqref{85}) that $\b^*>0$ and then $\a_2\neq 0$ follows as in \cite[Lemma 6.1]{BRR1}. In conclusion, it is apparent from \eqref{54} and \eqref{55} that, by taking $c$ sufficiently close to $c^*$, it is possible to take a component of the sets $\{\un{U}>0\}$ and $\{\un{V}>0\}$ in such a way \eqref{21} is satisfied, obtaining compactly supported generalized subsolutions satisfying \eqref{21}.
\end{proof}
\end{proposition}

\setcounter{equation}{0}
\setcounter{figure}{0}
\section{Asymptotic speed of propagation}
\label{sec6}
We are now able to give the proof of the first part of Theorem \ref{thmain}.

\vspace{0.2cm}
\noindent\emph{Proof of Theorem \ref{thmain}.} To prove the first condition of Definition \ref{deasp} we will use the supersolutions of \eqref{14} constructed in Proposition \ref{pr41} and Remark \ref{remark42}(i). We recall that they solve the linear system \eqref{41} and are of type
\begin{equation*}
(\bar{u}(x,t),\bar{v}(x,y,t))=e^{\a(x+c^*t)}(1,\g\phi(y))
\end{equation*}
with $\a,\g>0$ and $\phi(y)>0$ in $(0,L]$ satisfies $\phi'(0)>0$. As a consequence, since $(u_0,v_0)$ has compact support, there exists $k>0$ such that
\begin{equation*}
(u_0(x),v_0(x,y))<k(\bar{u}(x,0),\bar{v}(x,y,0))
\end{equation*}
for every $(x,y)\in\R\times(0,L]$. Moreover $(0,0)$ is a strict subsolution of \eqref{14}  and from Proposition \ref{prcp1} we have that
\begin{equation}
\label{61}
(0,0)\leq(u(x,t),v(x,y,t))<k(\bar{u}(x,t),\bar{v}(x,y,t))
\end{equation}
for every $t>0$. Observe that
\begin{equation*}
(\bar{\bar{u}}(x,t),\bar{\bar{v}}(x,y,t))=e^{\a(-x+c^*t)}(1,\g\phi(y))
\end{equation*}
is also a supersolution of \eqref{14} satisfying $(u_0(x),v_0(x,y))<k(\bar{\bar{u}}(x,0),\bar{\bar{v}}(x,y,0))$ for large $k$.

Fix now $c>c^*$, $t>0$, $|x|\geq ct$ and $y\in[0,L]$. We distinguish the cases $x\leq -ct<0$ and $-x\leq -ct<0$. In the first case, it follows that
\begin{equation*}
e^{\a(x+c^*t)}\leq e^{\a(c^*-c)t}
\end{equation*}
and this, together with \eqref{61}, implies
\begin{equation*}
(0,0)\leq(u(x,t),v(x,y,t))<k e^{\a(c^*-c)t}\left(1,\g\phi(y)\right)
\end{equation*}
for every $(x,y)\in\R\times(0,L]$ and \eqref{111} follows.  The second case can be treated  by comparing $(u,v)$ with $(\bar{\bar{u}},\bar{\bar{v}})$ in a similar fashion.

By adapting the arguments of the proof of Proposition \ref{pr34} to the case of Problem \eqref{51}, using the subsolution constructed in Proposition \ref{pr51}, it can be shown that
\begin{equation}
\label{62}
\lim_{t\to+\infty}(u(x+ct,t),v(x+ct,y,t))=\left(\frac{\nu}{\mu}V(L),V(y)\right)
\end{equation}
locally uniformly in $\overline{\O}$, where $V(y)$ is the unique solution of \eqref{110}. Property \eqref{112} now follows from \eqref{19} and \eqref{62} by using \cite[Lemma 4.4]{RTV}. \hfill\qedsymbol

\vspace{0.2cm}
Properties $(ii)$ and $(iii)$ of Theorem \ref{thmain} will be proved in the next sections.

\setcounter{equation}{0}
\setcounter{figure}{0}
\section{Limits for small and large diffusion on the road}
\label{sec7}

In this section we analyze the behavior of $c^*=c^*(D)$ as the diffusion on the road $D$ tends to $0$ and to $+\infty$, giving the proof of Theorem \ref{thmain}$(ii)$. The result regarding the first case is the following

\begin{proposition}
\label{pr71}
We have that
\begin{equation}
\label{71}
\lim_{D\da 0} c^*(D)=\ell_0>0,
\end{equation}
where $\ell_0$ is the asymptotic speed of propagation of the problem
\begin{equation}
\label{72}
  \left\{ \begin{array}{lll}
  u_t(x,t)=\nu v(x,L,t)-\m u(x,t) & \text{for } x\in\R, & t>0 \\
  v_t(x,y,t)-d\D v(x,y,t)=f(v) & \text{for } (x,y)\in\O, & t>0  \\ 
  d v_y(x,L,t)=\m u(x,t)-\nu v(x,L,t) & \text{for } x\in\R, & t>0 \\
  v(x,0,t)=0 & \text{for } x\in\R, & t>0. \end{array}\right.
\end{equation}

\begin{proof}
Observe first of all that limit \eqref{71} exists thanks to Theorem \ref{thmain}$(i)$. Fix now $D<d$. By the discussion of Sections \ref{sec4} and \ref{sec5}, we have that $c^*(D)$ is the one for which $\a_{D}^-(c,\b)$ given by \eqref{eqdef}
is tangent to $\a_{d}^+(c,\b)$ defined in \eqref{49}.
Moreover we have that necessarily the tangency occurs for $\b^*>\frac{\pi}{2L}$ (see Figure \ref{fig23}(B)). By passing to the limit for $D\da 0$ in \eqref{eqdef}, we get that $\ell_0$ is the unique value of $c$ for which $\a_{d}^+(c,\b)$ is tangent to
\begin{equation}
\label{73}
\a_0^-(c,\b)=\frac{-\mu d\b\cos(\b L)}{c\left(d\b\cos(\b L)+\nu\sin(\b L)\right)}, \quad \b\in\left[\frac{\pi}{2L},\ov{\b}\right),
\end{equation}
where $\ov{\b}$ is the one defined in \eqref{44}. The existence of such $c$ follows from the fact that
\begin{gather*}
\a_0^-\left(c,\frac{\pi}{2L}\right)=0, \qquad \a_0^-(c,\b)>0 \text{ for } \b\in\left(\frac{\pi}{2L},\ov{\b}\right), \qquad \lim_{\b\ua\ov{\b}}\a_0^-(c,\b)=+\infty, \\
\p_\b\a_0^-(c,\b)>0, \qquad \p_c\a_0^-(c,\b)<0, \\
\lim_{c\da 0}\a_0^-(c,\b)=+\infty, \qquad \lim_{c\ua +\infty}\a_0^-(c,\b)=0, 
\end{gather*}
together with the properties of $\a_d^+(c,\b)$ already described in Section \ref{sec4}.

To see that $\ell_0$ coincides with the asymptotic speed of propagation of \eqref{72}, it is sufficient, as in Section \ref{sec4}, to construct supersolutions of \eqref{72} of the form \eqref{42} for $c>\ell_0$ by intersecting the curves $\a_{d}^+(c,\b)$ and \eqref{73} and to proceed like in Section \ref{sec5} to construct compactly supported subsolutions to \eqref{72} for every $c<\ell_0$. Of course one has to prove the corresponding comparison principles for system \eqref{72}, which couples a strongly parabolic equation with a degenerate one. They essentially hold because, if there was a first contact point at a positive time between a supersolution and a subsolution, either it would be for the $v$ component, which is impossible since a classical comparison principle holds, or for the $u$ component at $y=L$. In such case, the time derivative of the difference between the super- and the subsolution would be negative, while the right-hand side of the first equation of \eqref{72} would be positive, obtaining again a contradiction (for a more detailed treatment of the comparison principles for such degenerate system in a similar context see \cite[Proposition 2.5]{RTV}).
\end{proof}
\end{proposition}

We now pass to the case $D\to+\infty$. 

\begin{proposition}
\label{pr72}
We have that $c^*(D)$ is unbounded as $D\to+\infty$ and
\begin{equation}
\label{74}
\lim_{D\ua \infty} \frac{c^*(D)}{\sqrt{D}}=\ell_\infty>0
\end{equation}
where $\ell_\infty$ is the asymptotic speed of propagation of the problem
\begin{equation}
\label{75}
  \left\{ \begin{array}{lll}
  u_t(x,t)-u_{xx}(x,t)=\nu v(x,L,t)-\m u(x,t) & \text{for } x\in\R,  & t>0 \\
  v_t(x,y,t)-d v_{yy}(x,y,t)=f(v) & \text{for } (x,y)\in\O,  & t>0  \\ 
  dv_y(x,L,t)=\m u(x,t)-\nu v(x,L,t) & \text{for } x\in\R, & t>0 \\
  v(x,0,t)=0 & \text{for } x\in\R, & t>0. \end{array}\right.
\end{equation}

\begin{proof}
Recall from the  proof of Proposition \ref{pr41} that, for large $D$,
\begin{equation}
\label{76}
c^*(D)=\min\{c^{*,1}(D),c^{*,2}(D)\},
\end{equation}
where $c^{*,1}$ is the first value of $c$ for which $\a_D^+(c,\b)$ and $\a_d^-(c,\b)$ intersect, being tangent and the same for $c^{*,2}$, considering $\tilde{\a}_D^+(c,\b)$ and $\tilde{\a}_d^-(c,\b)$. We will prove that \eqref{74} holds both for $c^{*,1}(D)$ and $c^{*,2}(D)$ and therefore \eqref{74} will follow from \eqref{76}.

We start with the case of $c^{*,1}(D)$ (for convenience, we will denote it by $c^{*}(D)$ when there is no possibility of confusion), which is increasing thanks to Proposition \ref{pr41} and admits a limit as $D\ua\infty$. It is obvious (see Figure \ref{fig26}(A)) that
\begin{equation}
\label{77}
\a_d^-(c^*(D),\tilde{\b}(c^*(D),D))<\a_D^+(c^*(D),0),
\end{equation}
where $\tilde{\b}(c^*(D),D)\in(\frac{\pi}{2L},\ov{\b})$ is the one defined in \eqref{45} (we have pointed out explicitly the dependence on $D$). Relation \eqref{77} can be written as
\begin{equation}
\label{78}
\frac{1}{d}\left(\!\!1\!-\!\sqrt{1-\frac{\cKPP^2-4d^2\tilde{\b}^2(c^*(D),D)}{{c^*}^2(D)}}\right)\!<\frac{1}{D}\left(1+\sqrt{1+\frac{4\mu d}{d+\nu L}\frac{D}{{c^*}^2(D)}}\right).
\end{equation}
Assume by contradiction that $c^*(D)$ is bounded. Then, from \eqref{45}, we have that
\begin{equation}
\label{79}
\frac{c^*(D)^2}{D}=\frac{-4\mu d\tilde{\b}(c^*(D),D)\cos(\tilde{\b}(c^*(D),D) L)}{d\tilde{\b}(c^*(D),D)\cos(\tilde{\b}(c^*(D),D) L)+\nu\sin(\tilde{\b}(c^*(D),D) L)},
\end{equation}
from which we get
\begin{equation*}
\lim_{D\ua\infty} \tilde{\b}(c^*(D),D)=\frac{\pi}{2L}.
\end{equation*}
By passing now to the limit as $D\ua\infty$ in \eqref{78}, we get a contradiction and, therefore,
\begin{equation}
\label{710}
\lim_{D\ua\infty}c^*(D)=+\infty.
\end{equation}
As the curves are tangent for the first time, we also have $\a_D^+(c^*,0)<\a_d^-(c^*,0)$, which reads
\begin{equation*}
\frac{1}{D}\left(1+\sqrt{1+\frac{4\mu d}{d+\nu L}\frac{D}{{c^*}^2}}\right)<\frac{1}{d}\left(1-\sqrt{1-\frac{\cKPP^2}{{c^*}^2}}\right).
\end{equation*}
Using \eqref{710}, we derive from the previous relation
\begin{equation}
\label{711}
\frac{1}{D}\left(1+\sqrt{1+\frac{4\mu d}{d+\nu L}\frac{D}{{c^*}^2}}\right)<{c^*}^{-2}\left(\frac{\cKPP^2}{2d}+o(1)\right),
\end{equation}
where, as usual, $o(1)$ denotes a quantity that goes to $0$ as $D\ua\infty$. Solving now for ${c^*}^2/D$, we obtain
\begin{equation*}
\frac{{c^*}^2}{D}<\frac{(d+\nu L)f'(0)^2}{(d+\nu L)f'(0)+\mu d}+o(1),
\end{equation*}
from which we conclude
\begin{equation}
\label{712}
\limsup_{D\ua\infty}\frac{{c^*}^2}{D}\leq\frac{(d+\nu L)f'(0)^2}{(d+\nu L)f'(0)+\mu d}.
\end{equation}
We now observe that
\begin{equation}
\label{713}
\tilde{\b}(c^*(D),D)<\sqrt{\frac{f'(0)}{d}}
\end{equation}
for every $D$, because otherwise
\begin{equation*}
\a_D^+(c^*(D),\tilde{\b}(c^*(D),D))>0\geq\a_d^-(c^*(D),\tilde{\b}(c^*(D),D))
\end{equation*}
and there would be another intersection between $\S_D(c^*)$ and $\S_d(c^*)$ apart from $\b^*$, which contradicts the construction of $c^*$.
As a consequence of \eqref{713} we have that
\begin{equation}
\label{714}
\liminf_{D\ua\infty}\tilde{\b}(c^*(D),D)\leq
\limsup_{D\ua\infty}\tilde{\b}(c^*(D),D)\leq\sqrt{\frac{f'(0)}{d}}.
\end{equation}
We now distinguish two cases. If
\begin{equation}
\label{715}
\liminf_{D\ua\infty}\tilde{\b}(c^*(D),D)=\sqrt{\frac{f'(0)}{d}},
\end{equation}
we have from \eqref{714} that
\begin{equation*}
\lim_{D\ua\infty}\tilde{\b}(c^*(D),D)=\sqrt{\frac{f'(0)}{d}}>\frac{\pi}{2L}
\end{equation*}
and, therefore, using \eqref{712} and taking the limsup for $D\ua\infty$ in \eqref{79}, we get that $\lim_{D\ua\infty}c^*(D)^2/D$
exists and is positive. On the other hand, if
\begin{equation}
\label{716}
\liminf_{D\ua\infty}\tilde{\b}(c^*(D),D)<\sqrt{\frac{f'(0)}{d}},
\end{equation}
we consider relation \eqref{78}, which, by using \eqref{710} and \eqref{716}, becomes
\begin{equation*}
2{c^*}^{-2}\left(f'(0)-d\tilde{\b}^2(c^*(D),D)+o(1)\right)<\frac{1}{D}\left(1+\sqrt{1+\frac{4\mu d}{d+\nu L}\frac{D}{{c^*}^2}}\right)
\end{equation*}
for large $D$. Solving for ${c^*}^2/D$, we now get
\begin{equation*}
\frac{{c^*}^2}{D}>\!\min\!\left\{\frac{(d+\nu L)(f'(0)-d\tilde{\b}^2(c^*(D),D))^2}{(d+\nu L)(f'(0)-d\tilde{\b}^2(c^*(D),D))+\mu d},2\left(f'(0)\!-\!d\tilde{\b}^2(c^*(D),D)\right)\right\}+o(1)
\end{equation*}
and, thanks to \eqref{716}, we have
\begin{equation}
\label{717}
\liminf_{D\to+\infty}\frac{{c^*}^2}{D}>0.
\end{equation}
Summing up, \eqref{717} holds both in case \eqref{715} and \eqref{716}.
This, together with \eqref{712}, implies that ${c^*}^2/D$ is bounded and bounded away from $0$. It is therefore natural to perform in \eqref{43} the change of variables
\begin{equation}
\label{718}
\hat{\a}=\sqrt{D}\a, \qquad \hat{c}=\frac{c}{\sqrt{D}},
\end{equation}
obtaining
\begin{equation}
\label{719}
\left\{\begin{array}{l}
-\hat{\a}^2+\hat{c}\hat{\a}+\frac{\mu}{1+\frac{\nu \tan(\b L)}{d \b}}=0 \\
-\frac{d\hat{\a}^2}{D}+\hat{c}\hat{\a}=f'(0)-d\b^2,
\end{array}\right.
\end{equation}
where $\hat{c}$ is bounded and bounded away from $0$. The first equation describes, in the plane $(\b,\hat{\a})$, the curve $\S_1(\hat{c})$ defined in \eqref{47}, therefore the function $\hat{\a}_1^+(\hat{c},\b)$, to which we will be interested in, is bounded and bounded away from $0$ for all $c>c_{\KPP}$ and $\b$ in the proper domain of definition. Therefore, by taking the limit for $D\ua\infty$ in \eqref{719}, we get
\begin{equation}
\label{720}
\left\{\begin{array}{l}
-\hat{\a}^2+\hat{c}\hat{\a}+\frac{\mu}{1+\frac{\nu \tan(\b L)}{d \b}}=0 \\
\hat{\a}=\frac{f'(0)-d\b^2}{\hat{c}}.
\end{array}\right.
\end{equation}
The second equation is a concave parabola, symmetric with respect to the $\hat{\a}-$axis, passes through $\left(\sqrt{f'(0)/d},0\right)$ and whose vertex is $\left(0,f'(0)/\hat{c}\right)$.

Now we pass to the case of $c^{*,2}$, which, as above, will be simply denoted by $c^*$. In this case, as it is apparent from Figure \ref{fig26}(B), we have
\begin{equation*}
\tilde{\a}_d^-(c^*(D),0)<\lim_{\b\to+\infty}\tilde{\a}_D^+(c^*(D),\b),
\end{equation*}
which reads as
\begin{equation}
\label{721}
\frac{1}{d}\left(1-\sqrt{1-\frac{\cKPP^2}{{c^*}^2}}\right)<\frac{1}{D}\left(1+\sqrt{1+4\mu \frac{D}{{c^*}^2}}\right)
\end{equation}
and gives that
\begin{equation*}
\lim_{D\to+\infty}c^*(D)=+\infty,
\end{equation*}
because, if the limit was finite, say $\ell$, passing to the limit in \eqref{721} would lead to
\begin{equation*}
0<\frac{1}{d}\left(1-\sqrt{1-\frac{\cKPP^2}{\ell^2}}\right)\leq 0,
\end{equation*}
which is impossible. Therefore, \eqref{721} gives
\begin{equation*}
{c^*}^{-2}\left(\frac{\cKPP^2}{2d}+o(1)\right)<\frac{1}{D}\left(1+\sqrt{1+4\mu \frac{D}{{c^*}^2}}\right)
\end{equation*}
and, solving for ${c^*}^2/D$ and taking the liminf as $D\ua\infty$ we get
\begin{equation*}
\liminf_{D\to+\infty}\frac{{c^*}^2}{D}\geq\min\left\{\frac{f'(0)^2}{f'(0)+\mu},2f'(0)\right\}.
\end{equation*}
On the other hand, in this situation we also have that $\tilde{\a}_D^+(c^*(D),0)<\tilde{\a}_d^-(c^*(D),0)$,
which provides us with \eqref{711} and, therefore, with the upper bound for ${c^*}^2/D$ given by \eqref{712}.
By performing the change of variables \eqref{718} in \eqref{413} and passing to the limit for $D\to+\infty$ we obtain
\begin{equation}
\label{722}
\left\{\begin{array}{l}
-\hat{\a}^2+\hat{c}\hat{\a}+\frac{\mu}{1+\frac{\nu \tanh(\b L)}{d \b}}=0 \\
\hat{\a}=\frac{f'(0)+d\b^2}{\hat{c}}.
\end{array}\right.
\end{equation}
The first equation describes, in the plane $(\b,\hat{\a})$, the curve $\tilde{\S}_1(\hat{c})$ defined in \eqref{418}, while the second one is a parabola, which is symmetric to the one of the second equation of \eqref{720} with respect to the line $\hat{\a}=f'(0)/\hat{c}$.

With a similar reasoning as that of Section \ref{sec4}, it is easily seen that there is a smallest value of $\hat{c}$ for which either the two curves of \eqref{720} or of \eqref{722} are tangent, which provides us with $\ell_\infty$.

To see that this limit coincides with the asymptotic speed of propagation of Problem \eqref{75} it suffices to repeat the construction of Sections \ref{sec4}-\ref{sec5} for this problem, starting from supersolutions of type \eqref{42} and \eqref{412} and using comparison principles analogous to the ones of Section \ref{sec2}, which can be proved for \eqref{75} by applying the parabolic maximum principle in $y$ on every slice, with fixed $x$ (see \cite{RTV} for a detailed proof in a similar context and \cite{D}, in the context of travelling waves, for comparison principles related to this degenerate system with Neumann boundary conditions at $y=0$).
\end{proof}
\end{proposition}

\setcounter{equation}{0}
\setcounter{figure}{0}
\section{Influence of the road and limit for large field}
\label{sec8}

To examine the influence of the road in Problem \eqref{14}, it is appropriate to compare its asymptotic speed of propagation with the one of the following problem
\begin{equation}
\label{81}
  \left\{ \begin{array}{lll}
  v_t(x,y,t)-d\D v(x,y,t)=f(v) & \text{for } (x,y)\in\O, & t>0  \\ 
  dv_y(x,L,t)=-\nu v(x,L,t) & \text{for } x\in\R, & t>0 \\
  v(x,0,t)=0 & \text{for } x\in\R, & t>0, \end{array}\right.
\end{equation}
which models a classical Fisher-KPP diffusion in the strip $\O$ and part of the population $v$ just leaves the field at level $y=L$.

By using the same techniques of Section \ref{sec3} it is possible to show that Problem \eqref{81} admits a unique positive steady state if and only if
\begin{equation}
\label{82}
\frac{f'(0)}{d}>\ov{\b}^2>\left(\frac{\pi}{2L}\right)^2,
\end{equation}
where $\ov{\b}=\ov{\b}(d,\nu,L)$ is, as in Section \ref{sec4}, the first positive value for which \eqref{44} vanishes. By comparing \eqref{82} with \eqref{17}, it is apparent that one effect of the road is to enhance the persistence of the species, since the condition for persistence is less restrictive in the presence of the road. On the other hand, when \eqref{82} holds, by taking
\begin{equation*}
\ov{v}(x,y,t)=e^{\a(x+ct)}\sin(\ov{\b} y)
\end{equation*}
as supersolution and following the lines of Sections \ref{sec4}--\ref{sec6}, it is possible to show that Problem \eqref{81} admits an asymptotic speed of propagation
\begin{equation*}
c_{\KPP}^{\DR}=2\sqrt{d\left(f'(0)-d\ov{\b}^2\right)}<c_{\KPP}
\end{equation*}
(here $\DR$ stands for the Dirichlet-Robin boundary conditions associated to the Fisher-KPP equation in \eqref{81}).

Recalling the monotonicity property of $c^*(D)$ given by Theorem \ref{thmain}$(i)$, we have that $c^*(D)>\ell_0$, where $\ell_0$ is the one constructed in Proposition \ref{pr71} as the smallest value of $c$ for which the curves \eqref{73} and $\a_d^+(c,\b)$ intersect. Observe that the former is defined for $\b<\ov{\b}$, while, recalling \eqref{410}, $\a_d^+(c_{\KPP}^{\DR},\b)$ is defined for $\b\geq\ov{\b}$. This means that, for every $D$,
\begin{equation}
\label{83}
c^*(D)\geq\lim_{D\da 0}c^*(D)>c_{\KPP}^{\DR}.
\end{equation}
Therefore, a second effect of the road is speeding up the propagation in the field. Moreover, from the second relation of Theorem \ref{thmain}$(ii)$, we have that this effect can be arbitrarily enhanced, provided that $D$ is sufficiently large.

We conclude with the proof of Theorem \ref{thmain}, considering the limit for large field. We will emphasize the dependence of $c^*$ on the width of the strip $L$, by writing $c^*_L$. 

\vspace{0.2cm}
\noindent\emph{Proof of Theorem \ref{thmain}(iii).} We recall that in \cite{BRR1} it was proved that
 \begin{equation*}
c^*_{\infty}(D)\begin{cases}
=\cKPP &\text{if } D\leq 2d \\
>\cKPP &\text{if } D> 2d.
\end{cases}
\end{equation*}
We distinguish the same two cases, starting with $D\leq 2d$. From \eqref{83} and the discussion of Section \ref{sec4} we have that
\begin{equation}
\label{84}
c_{\KPP}^{\DR}(L)< c^*_L<\cKPP.
\end{equation}
Recalling from \eqref{44} that  $\ov{\b}(L)\in\left(\frac{\pi}{2L},\frac{\pi}{L}\right)$, we have that \eqref{84} implies
\begin{equation*}
\lim_{L\to+\infty}c^*_L=c_{\KPP}.
\end{equation*}

We now assume $D>2d$ and recall (see \cite{BRR1}) that $c^*_{\infty}$ is the value of $c$ for which the curve
\begin{equation*}
\tilde{\a}_D^{\infty}(c,\b):=\frac{1}{2D}\left(c+\sqrt{c^2+\frac{4\mu D}{1+\frac{\nu}{d\b}}}\right),
\end{equation*}
and the curve $\tilde{\a}_d^-(c,\b)$ defined in \eqref{419} are tangent.

In our case, a crucial role will be played by the behavior of the curves $\S_D(c)$, $\S_d(c)$, $\tilde\S_D(c)$, $\tilde\S_d(c)$ for $c=\cint$, where $\cint$ is the one of \eqref{420}. It is possible to show with direct computations that $D\mapsto \cint(D)$ is increasing for $D>2d$. As a consequence, there exists $\tilde D$ such that $\cint(D)>\cKPP$ for every $D>\tilde D$. Due to the monotonicity of $\cint(D)$, we deduce that $\tilde D$ is the unique value of $D$ for which $\a_D^+(\cKPP,0)=\a_d^\pm(\cKPP,0)$, which, taking into account \eqref{eqstar} and \eqref{ipdeg}, provides
\begin{equation*}
\tilde D=\tilde D(L)=2d+\frac{4\mu d^3}{(d+\nu L)\cKPP^2}.
\end{equation*}
In particular, since we are taking $D>2d$, we have that, for $L$ large enough, $D>\tilde D(L)$ and, therefore, $\cint>\cKPP$.
Moreover, recalling \eqref{tipo2} and \eqref{419},
\begin{equation}
\label{85}
\p_{\b\b}\left(\tilde{\a}_{D}^+(\cint,0)-\tilde{\a}_{d}^-(\cint,0)\right)=2d\left(\frac{\psi}{3\zeta\sqrt{\cint^2+\frac{4\mu d D}{d+\nu L}}}-\frac{1}{\sqrt{\cint^2-\cKPP^2}}\right),
\end{equation}
where we have set $\psi:=L^3\mu\nu$ and $\zeta:=(d+\nu L)^2$. It is possible to check, by direct computations, that \eqref{85} has the same sign as
\begin{multline*}
(\psi^2-9\zeta^2) \left((d+\nu L)D\cKPP^2+ 4\mu d^3\right)^2+\\
-4d(D-d)\left((d+\nu L)\cKPP^2 + 4\mu d^2\right)\left(36\zeta^2\mu d D+\psi^2(d+\nu L)\cKPP^2\right),
\end{multline*}
which is positive for large $L$, since it is a polynomial of degree $8$ in $L$ with leading coefficient equal to $\left(\mu\nu^2c_{\KPP}^2(D-2d)\right)^2>0$. This implies that
\begin{equation*}
\tilde{\a}_{D}^+(\cint,\b)>\tilde{\a}_{d}^-(\cint,\b) \qquad \text{ for } \b>0, \quad \b\sim0
\end{equation*}
and, since the curve $\tilde\S_d(\cint)$ intersects the $\a-$axis in another point, which lies above the intersection with $\tilde\S_D(\cint)$, we obtain that $\tilde{\S}_D(\cint)$ and $\tilde{\S}_d(\cint)$ intersect, apart from $\b=0$, in a positive value of $\b$ too.

On the other hand, using \eqref{420} and defining
\begin{equation*}
\mathfrak{c}^2:=\lim_{L\to+\infty}\cint^2(L)=\frac{D^2c_{\KPP}^2}{4d(D-d)},
\end{equation*}
it is easy to see that the curve $\S_D(c^*_L)$ introduced in \eqref{47} approaches in the $(\b,\a)-$plane, as $L$ goes to $+\infty$, the vertical segment $\{0\}\times [0,\mathfrak{c}/D)$. As a consequence of these considerations and the discussion of Section \ref{sec4}, we have that, for large $L$, $c^*_L$ is obtained as the value for which $\tilde{\a}_D^+(c,\b)$ and $\tilde{\a}_d^-(c,\b)$ are tangent. Observe that $\tilde{\a}_D^+$ is decreasing in $L$ and, therefore, $c^*_L$ is increasing. In addition $\tilde{\a}_D^+>\tilde{\a}_D^{\infty}$,
which entails that $c^*_L<c^*_\infty$. Finally, $\tilde{\a}_D^+$ tends, as $L\to+\infty$, to $\tilde{\a}_D^{\infty}$, together with its derivatives, locally uniformly in $\R_+\times\R_+$. This implies that
\begin{equation*}
\lim_{L\to+\infty}c^*_L=c^*_\infty
\end{equation*}
also in this case and concludes the proof. \hfill \qedsymbol

\begin{remark}
At a first glance, one may think that the result provided in Theorem \ref{thmain}(iii) gives a contradiction for $D\in(2d,\limsup_{L\to+\infty}D_{\KPP})$, since \eqref{84} holds for $D<D_{\KPP}$ and one could repeat the argument of the first part of the proof to show that $c^*_L$ converges to $\cKPP$ as $L\to+\infty$. However, this is not the case, as the previous interval is empty, since
\begin{equation}
\label{86}
\lim_{L\to+\infty}D_{\KPP}(L)=2d.
\end{equation}
This follows as a byproduct of the proof of Theorem \ref{thmain}(iii), but for clarity we provide a direct proof here. Thanks to the monotonicity in $L$ and $D$ of $\S_D(\cKPP)$ and by definition of $D_{\KPP}$, it follows easily that $L\mapsto D_{\KPP}(L)$ is decreasing and therefore the limit in \eqref{86} exists. Let us denote it by $\ov D$ and assume by contradiction that it satisfies $\ov D>2d$. As shown in the last part of the proof of Theorem \ref{thmain}(iii), the curve $\S_{D_{\KPP}(L)}(\cKPP)$ converges, as $L\to+\infty$, to the segment $\{0\}\times [0,\cKPP/\ov D)$, which, since we are assuming $\ov D>2d$, lies at a positive distance from $\S_d(\cKPP)$. As a consequence, there would be no intersection between the latter curve and $\S_{D_{\KPP}(L)}(\cKPP)$ for large $L$ either, contradicting the definition of $D_{\KPP}$.
\end{remark}

\section*{Acknowledgement}
The author is deeply grateful to Professors Henri Berestycki, Jean-Michel Roquejoffre and Luca Rossi for the broad discussions that helped  him in the preparation of this work. Moreover he wishes to thank all the people at CAMS (\emph{Centre d'analyse et de math\'ematique sociales}, Paris), where most of the work was developed, in 2013, for the very warm treatment.

\end{document}